\documentclass{amsart}
\usepackage[utf8]{inputenc}

\usepackage{hyperref}
\usepackage{enumitem}
\usepackage{tikz-cd}
\tikzset{
symbol/.style={
draw=none,
every to/.append style={
edge node={node [sloped, allow upside down, auto=false]{$#1$}}}
}
}
\usepackage{pgfplots}
\pgfplotsset{compat=newest}
\usepackage{mathtools,amsmath,amssymb}
\usepackage{bbm}

\usepackage{tikz}
\usetikzlibrary{cd,arrows}
\tikzcdset{scale cd/.style={every label/.append style={scale=#1},
cells={nodes={scale=#1}}}}
\usetikzlibrary{decorations.markings}

\usepackage[backend=biber, style=alphabetic, giveninits=true]{biblatex} 
\theoremstyle{plain}
\newtheorem{thm}{Theorem}[section]
\newtheorem{lemma}[thm]{Lemma}
\newtheorem{prop}[thm]{Proposition}
\newtheorem{cor}[thm]{Corollary}

\theoremstyle{definition}
\newtheorem{defi}[thm]{Definition}

\newtheorem{rmk}[thm]{Remark}

\newtheorem{nota}[thm]{Notation}

\newcommand{\N}{\mathbb{N}}
\newcommand{\Z}{\mathbb{Z}}
\newcommand{\R}{\mathbb{R}}
\newcommand{\1}{\mathbbm{1}}
\newcommand{\Oe}{{\Omega_\e}}
\newcommand{\Oet}{{\Omega_\e(t)}}
\newcommand{\e}{\varepsilon}
\newcommand{\intS}{\int\limits_S}
\newcommand{\intO}{\int\limits_\Omega}
\newcommand{\intOe}{\int\limits_\Oe}

\newcommand{\intY}{\int\limits_Y}
\newcommand{\intYp}{\int\limits_{\Yp}}
\newcommand{\intYpx}{\int\limits_{\Ypx(t)}}
\newcommand{\intOY}{\intO \!\! \intY}
\newcommand{\intOYp}{\intO \!\intYp}
\newcommand{\intOYpx}{\intO \! \intYpx}
\newcommand{\intSOYp}{\intSO \! \intYp}
\newcommand{\intG}{\int\limits_\Gamma}
\newcommand{\intSO}{\intS \!\! \intO}
\newcommand{\intSOe}{\intS \!\! \intOe}
\newcommand{\intSOY}{\intS \!\! \intOY}
\newcommand{\norm}[2]{\left|\left|#1\right| \right|_{#2}}
\newcommand{\Yp}{{Y^*}}
\newcommand{\Ypx}{{Y^*_x}}
\newcommand{\cYp}{\chi_\Yp}
\newcommand{\hu}{{\hat{u}}}
\newcommand{\hue}{{{\hat{u}_\e}}}

\newlength{\myl}
\newcommand{\TscW}[1]{
\settowidth{\myl}{$\xrightharpoonup{\hspace{0.05cm} #1 \hspace{0.05cm}}$}
\xrightharpoonup{\hspace{0.05cm} #1 \hspace{0.05cm}}\hspace{-\myl} \xrightharpoonup{\hspace{0.05cm}\phantom{#1} \hspace{0.15cm}}
}

\newcommand{\TscS}[1]{
\settowidth{\myl}{$\rightarrow{ #1 }$}
\xrightarrow{ #1 \hspace{0.05cm}}
\hspace{-1.03\myl} \ \ 
\xrightarrow{\hspace{0.05cm}\phantom{#1} \hspace{0.15cm}}
}

\newcommand{\Ge}{{\Gamma_\e}}
\newcommand{\Te}{{\mathcal{T}_\e}}
\renewcommand{\div}{{\operatorname{div}}}
\newcommand{\tx}{(t,x)}
\newcommand{\txy}{(t,x,y)}
\newcommand{\Gek}{{\Gamma_{\e,k}}}
\newcommand{\Gekt}{{\Gamma_{\e,k}(t)}}
\newcommand{\limEps}{\lim\limits_{\e \to 0}}
\newcommand{\txxeps}{\left(t,x, \tfrac{x}{\e}\right)}
\newcommand{\SO}{{S\times \Omega}}
\newcommand{\xYp}{{\chi_{\Yp}}}

\newcommand{\xYpt}{\chi_{Y^*\psi_0^{-1}}}
\newcommand{\dx}{{dx}}
\newcommand{\dt}{{dt}}
\newcommand{\dy}{{dy}}
\newcommand{\dxt}{{\dx \dt}}
\newcommand{\dyxt}{{\dy \dx \dt}}
\newcommand{\p}{{\mathrm{p}}}
\newcommand{\dd}{{d}}

\newcommand{\A}{{\mathcal{A}}}
\newcommand{\B}{{\mathcal{B}}}
\newcommand{\V}{{\mathcal{V}}}
\newcommand{\W}{{\mathcal{W}}}

\newcommand{\PorosityA}[3]{
({0.23 + 0.12*sin((#3+ 0.3*#1 + 0.2*#2) r)})
}

\newcommand{\Ol}{22}
\newcommand{\Ow}{12}

\newcommand{\OetGrafik}[1]{
\OmegaScope{
\foreach \x in {-1,...,\Ol} {
\foreach \y in {-1,...,\Ow} {
\Cell{\x}{\y}{\PorosityA{\x}{\y}{#1}}
}}}
}

\newcommand{\Cell}[3]{
\draw [fill = gray] (#1 +0.5, #2 + 0.5) circle (#3);
}

\newcommand{\CellWithBoundary}[3]{
\draw	(#1, #2) rectangle (#1 +1, #2 +1);
\Cell{#1}{#2}{#3};
}

\newcommand{\OmegaScope}[1]{
\begin{scope}
\path[draw, clip] (2,0) -- (6,0) --(6,1) --(10,1) --(10,0)-- (12,0) --(12,2) --(14,2) -- (14,6) -- (15,6) --(15,10) --(14,10) -- (14,12) -- (12,12)--(12,11) -- (10,11) -- (8,11)  -- (8,12) -- (4,12) --(4,10) --(0,10) --(0,2) --(2,2) --(2,0);
{#1}
\end{scope}
}

\title[Homogenisation of local colloid evolution]{Homogenisation of local colloid evolution induced by reaction and diffusion}
\author{David Wiedemann}
\author{Malte A. Peter}
\thanks{D.W. was partially supported by a doctoral scholarship
provided by the Studienstiftung des deutschen Volkes.}
\address{(D.W.) Institute of Mathematics, University of Augsburg, 86135 Augsburg, Germany}
\email{david.wiedemann@math.uni-augsburg.de}
\address{(M.A.P.) Institute of Mathematics, University of Augsburg, 86135 Augsburg, Germany and Centre for Advanced Analytics and Predictive Sciences (CAAPS), University of Augsburg, 86135 Augsburg, Germany}
\email{malte.peter@math.uni-augsburg.de}


\begin{document}
\begin{abstract}
We consider the homogenisation of a coupled reaction--diffusion process in a porous medium with evolving microstructure. A concentration-dependent reaction rate at the interface of the pores with the solid matrix induces a concentration-dependent evolution of the domain. Hence, the evolution is fully coupled with the reaction--diffusion process. In order to pass to the homogenisation limit, we employ the  two-scale-transformation method. Thus, we homogenise a highly non-linear problem in a periodic and in time cylindrical domain instead.
The homogenisation result is a reaction--diffusion equation, which is coupled with an internal variable, representing the local evolution of the pore structure.	
\end{abstract}
\maketitle
\tableofcontents
\section{Introduction}
Reaction--diffusion mechanisms in porous media often induce an evolution of the solid matrix. Typical examples are reaction mechanisms producing or consuming constituents which are part of the solid matrix, e.g. in concrete carbonation (cf.~\cite{Kro95}, \cite{Bie88}) or crystal precipitation and dissolution (cf.~\cite{TSR07}, \cite{VNo08}). Similarly, if biofilms are present, these can often be viewed as a solid-matrix-type part of the porous medium on the pore scale. In this context, production of biofilm can be modelled on the microscale similarly as production of solid matrix (cf.~\cite{TZK02}, \cite{KAP07}, \cite{VNP10}).

Mathematical models for reaction and diffusion in porous media are typically obtained from upscaling processes on the pore scale by averaging or homogenisation techniques. A classic method in this context is periodic homogenisation (cf.~\cite{All92}, \cite{Ngu89}), which has been extended to cope with (non-periodic) evolving microstructures (cf.~\cite{Pet07}). The extension relies on transforming the non-periodic evolution to a periodic reference geometry, which requires modelling of this (concentration-dependent) transformation in the context of particular applications, for instance a detailed discussion for concrete carbonation can be found in \cite{PB09a}.
\begin{equation}\label{diag}
\hspace{-0.1cm}
\begin{tikzcd}[scale cd=0.82, row sep=1cm, column sep = 4.9cm]
\textrm{evolving microproblem} \arrow[r, "\textrm{homogenisation on the evolving domain}"] \arrow[d, "\textrm{transformation}" ] & \textrm{evolving macroproblem} 
\arrow[d, leftarrow, "\textrm{back-transformation}"] \\
\textrm{transformed microproblem}
\arrow[r, "\textrm{homogenisation on periodic reference domain}"]& \textrm{transformed macroproblem}
\end{tikzcd}
\end{equation}
The approach of transforming on a periodic reference domain has found also application in the homogenisation of thermoelasticity \cite{EM17} or the homogenisation of advection--reaction--diffiusion problems in porous media (cf.~\cite{GNP21}), where the domain’s evolution is a priori given.
Moreover, it has been recently shown that the homogenisation of the substitute problem is equivalent to the homogenisation of the actual problem in the non-periodic mirostructure, i.e.~that \eqref{diag} commutes  (cf.~\cite{Wie21}). Furthermore, a new two-scale-transformation rule has been been derived there, which yields a transformation-independent homogenisation result after the back-transformation. 

In the present paper, we use this approach to homogenise rigorously a reaction--diffusion problem where the domain evolution is not a priori given but coupled with the solution itself. 
The homogenisation of problems where the evolving microstructure is coupled with the solution itself has been also considered by a level-set approach. There, the domain is described by a level-set function solving a level-set equation, which involves the other unknowns. In this framework, microscopic models for crystal precipitation and dissolution (cf.~\cite{VNP10}) or biofilm growth in porous media (cf.~\cite{SRF16}) have been homogenised. However, the corresponding effective macroscopic problems have been derived by formal asymptotic expansion only. Numerical simulations and analytical discussion of such type of limit models can be found in \cite{GFP+20}, \cite{GFK+22}, \cite{KGB+22}.

In this manuscript, we revisit the microscale model by \cite{MN20} for one reaction--diffusion equation and derive their upscaled model by a mathematically rigorous homogenisation procedure based on the recent results of \cite{Wie21}. In this context, we show that such coupling of the pore structure with the solution of the reaction--diffusion equation can be handled by the two-scale-transformation method. For this purpose, we construct a concrete $\e$-scaled transformation for the $\e$-scaled domains by means of a generic parametrisable cell transformation. There, the radius of the solid obstacles becomes the parameter. By showing a certain kind of strong convergence for the radii of the $\e$-scaled model, we can verify the assumptions of the two-scale-transformation method. Thus, we can pass rigorously to the two-scale limit in the substitute problem. Moreover, using the two-scale-transformation rule of \cite{Wie21}, we obtain a two-scale limit problem in the actual non-cylindrical evolving two-scale domain, which is independent of the chosen transformation. There, we split the macroscopic and microscopic variables in order to derive an effective equation. The result is a macroscopic reaction--diffusion problem coupled with an internal variable, which represents the local radius of the solid. This local radius is given by an ordinary differential equation and scales not only the time-derivative term and the reaction rate of the reaction--diffusion equation but also affects the effective diffusivity. The diffusivity is still computed by solutions of cell problems as in the case of a rigid domain. However, the domain for the cell problems is now parametrised by the internal radius and affects in this way the local effective diffusivity.

This paper is organised as follows: In section 2, we derive the microscopic model \eqref{eq:StrongForm1}--\eqref{eq:StrongForm4}, which consists of a reaction--diffusion problem coupled with the evolution of the domain. Then, we state the corresponding weak formulation in the evolving domain. Using a generic cell transformation, we transform the weak form to the equivalent weak form \eqref{eq:def:psi-eps}--\eqref{eq:def:Psi-eps}, \eqref{eq:WeakTrans1}--\eqref{eq:WeakTrans3} on the periodic substitute domain, which becomes highly non-linear.
In section 3, we show the existence and uniqueness of the solution of the transformed microscopic model by a fixed point argument. There, we utilise the assumption that the radii, which define the solid domain, are a priori bounded from below and above. Moreover, we derive some $\e$-independent a priori estimates.
In section 4, we use two-scale convergence in order to pass to the homogenisation limit.
Since the coefficients in the equation depend on the solution itself, the problem becomes highly non-linear and we need a strong convergence of the solution. However, we can not not derive easily a uniform bound of the time derivative of the solution of the diffusion equation. Therefore, we can not use the classical Aubin--Lions lemma. Instead, we shift the solution of the reaction--diffusion equation with respect to time and estimate the difference to the actual solution. Then, we can conclude with the Simon-Kolmogorov compactness criterion (cf.~\cite[Theorem 1]{Sim86}) the strong convergence of the concentration.
Using this strong convergence, we can show a strong convergence of the radii, which allows us to apply the two-scale-transformation method. Thus, we can derive the two-scale limit problem in the cylindrical two-scale reference domain rigorously. In section 5, we transform the limit problem back and obtain the transformation-independent two-scale limit problem.
Then, we split the macroscopic and microscopic variable. This gives the effective problem
\eqref{eq:WeakTwoScaleLimit2}, \eqref{eq:homogenised-eq}--\eqref{eq:CellProblem} with its cell problems, which depend on an internal variable representing the local radius.

We use the following notations. Let $f,g \in L^2(U)$ and $U \subset \R^m$ for $m \in \N$, then we write the scalar product and the norm by:
$
(f,g)_U \coloneqq \int_U f(x)g(x) \dx$, $\norm{U}{U}^2 \coloneqq (f,f)_U$.
For $f \in H^1(U)'$ and $g \in H^1(U)$, we write the dual paring by $\langle f, g \rangle_\Omega \coloneqq \langle f, g \rangle_{H^1(U)',H^1(U)}$.

Furthermore, we use $C$ as generic constant which is independent of $\e$ and other variables and depends only on fixed constants. In cases, in which the generic constant can depend on other variables as for instance $\e$, we mark this by a subscript, e.g.~we write $C_\e$.
Moreover, let the spatial dimension be $N \in \N$ with $N \geq 2$.
\section{The mathematical model}
Let $\Omega$ be an open set in $\R^N$, which represents the macroscopic domain of the porous medium and let $\e = (\e_n)_{n \in \N}$ be a positive sequence converging to zero. We assume that $\Omega$  consists of whole $\e$-scaled cells $Y \coloneqq (0,1)^N$, i.e.~$\Omega = \operatorname{int}\left( \bigcup_{k \in I_\e} \e k +  \overline{Y}\right)$ for $I_\e \coloneqq \{k \in \Z^N \mid \e k + \e Y \cap \Omega \neq \emptyset\}$.
Moreover, we assume in the following that $\e\leq 1$.

\begin{figure}[h]
\centering
\begin{tikzpicture}[scale = 0.25]
\OetGrafik{0.2}

\begin{scope}[shift={(20,0)}]
\OetGrafik{4}
\end{scope}
\end{tikzpicture}
\caption{\label{fig:Oet_t=0_>0}The domain $\Oe(t)$ for $t = 0$ (left) and $t >0$ (right)}
\end{figure}
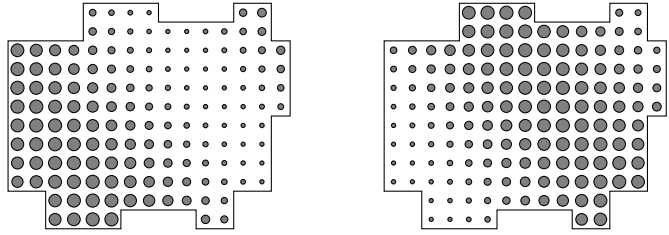

We assume that the pore structure of the porous medium is given by spherical obstacles in the cells $\e k + \e Y $ for $k \in I_\e$ which can grow and shrink on the time interval $S = (0,T)$ with $0<T<\infty$.
Thus, the $\e$-scaled porous medium is defined by 
\begin{align}\label{eq:def:Oet}
\Oe(t) \coloneqq \Omega \setminus\bigcup\limits_{k \in I_\e} \e \overline{B_{r_{\e,k}(t)}(k + x_M)}
\end{align}
where $x_M \coloneqq (0.5, \dots, 0.5)^\top$ is the centre of the reference cell and $r_{\e,k}(t)$ is the $\e^{-1}$-scaled radius of the solid obstacle in the cell located at $\e k$ at time $t \in S$ (cf.~Figure \ref{fig:Oet_t=0_>0}).

We assume that the size of the obstacles $ \e B_{r_{\e,k}(t)}(k+x_M) $ is affected by reactions on their surfaces $\Gekt \coloneqq \partial  \e B_{r_{\e,k}(t)}(k+x_M) $. The reactions rate $\e f(u_\e\tx,r_{\e,k}(t))$ depends on the concentration rate $u_\e$ and on the radius $r_{\e,k}$ of $S_{\e,k}$. Because the reaction rate depends on the radius, we can ensure $r_\mathrm{min}\leq r_{\e,k}(t) \leq  r_\mathrm{max}$ for every $k \in I_\e$ and every $t \in S$ for constants $0< r_\mathrm{min}< r_\mathrm{max}<0.5$ by the assumptions:
\begin{align}\label{eq:f>0}
&f(\cdot, r) \geq 0 \textrm{ for } r \leq r_\mathrm{min},\\\label{eq:f<0} &f(\cdot, r) \leq 0 \textrm{ for } r \geq r_\mathrm{max}.
\end{align}
Moreover, we assume that $f$ is uniformly Lipschitz continuous and bounded, i.e.~there exists a constant $C$ such that
\begin{align}\label{eq:f-Lipschitz}
&f(u_2,r_2) - f(u_2,r_1) \leq C (|u_2- u_1| +|r_2-r_1|),
\\\label{eq:f-bounded}
&|f(u_1,r_1)| \leq C_f
\end{align}
for $u_1, u_2  \in \R$ and $r_1,r_2 \in \R$.

We consider the case that the formed or vanishing solid has a constant concentration density $c_s$. Thus, the conservation of mass yields
\begin{align}\label{eq:MassConservation}
\tfrac{d}{dt} |\e B_{r_ {\e,k}(t)}(k+ x_M) | c_s
= \int\limits_{\Gekt} j_\e\tx \cdot n\tx d \sigma_x
\end{align}
where $j_\e\tx$ is the flux through $\Gekt$ and $n$ is the inner unit normal of $\e B_{r_{\e,k}(t)}(k+x_M)$.
We note that this flux consists of the diffusive flux and a flux which is induced by the evolution of the domain. 
We model the diffusive flux $j_D = -D \nabla u_\e\tx$ by Fick's law with a diffusion coefficient $D$.
The second flux, which is induced by the evolution of the domain, can be understood in the following sense: when the carrier medium becomes solid any excess dissolved concentration separates from the carrier medium and is pushed away, i.e.~$j_\Gamma\tx =  - v_{\Gamma_{\e,k}\tx} u\tx$, where $v_{\Gamma_{\e,k}}$ is the velocity of the boundary deformation. We note that $v_{\Gamma_{\e,k}}$ can be formulated explicitly by $v_{\Gamma_{\e,k}}\tx  = -\e \partial_t r_{\e,k}(t)  n\tx$.
Thus, the total flux on the boundary is
\begin{align}\label{eq:Derivation:boundCond}
j\tx = j_D\tx+ j_\Gamma\tx = -D \nabla u_\e\tx - v_{\Gamma_{\e,k}}\tx u\tx 
\end{align}
for $ t \in S$ and $ x \in \Gekt$.
On the other hand, the flux at the boundary in the normal direction, $j\tx \cdot n\tx$, represent the consumption or gain of concentration due to the reactions on $\Gekt$, which yields
\begin{align}\label{eq:BoundaryCondition}
(-D \nabla u_\e\tx - v_{\Gamma_{\e,k}}\tx u\tx )\cdot n\tx = j\tx \cdot n\tx = \e f(u_\e\tx,r_{\e,k}(t))
\end{align}
and equivalently
\begin{align}\label{eq:BoundaryCondition2}
-D \nabla u_\e\tx\cdot n\tx + \e \partial_t r_{\e,k}(t) u\tx  = \e f(u_\e\tx,r_{\e,k}(t)).
\end{align}
Inserting \eqref{eq:BoundaryCondition} in \eqref{eq:MassConservation} yields
\begin{align}\label{eq:d_tS(t)}
\tfrac{d}{dt} |\e B_{r_ {\e,k}(t)}(k+ x_M)| c_s = \int\limits_{\Gekt} \e f(u_\e\tx,r_{\e,k}(t)) d \sigma_x
\end{align}
and elementary calculus implies
\begin{align*}
\tfrac{d}{dt} |\e B_{r_ {\e,k}(t)}(k+ x_M)| = \e^N \tfrac{d}{\dt} V_N ( r_{\e,k}(t)) 
= 
\e^N S_{N-1} (r_{\e,k}(t)) \partial_t r_{\e,k}(t),
\end{align*}
where $V_N(r)$ denotes the volume of the $N$-ball with radius $r$ and $S_N(r)$ denotes the surface of the $N$-sphere with radius $r$.
Thus, we obtain the following ordinary differential equation for the radii:
\begin{align}\label{eq:Derivation:Partial:r}
\partial_t r_{\e,k}(t) =  \tfrac{\e^{-N} }{c_sS_{N-1}(r_{\e,k}(t)) } \int\limits_{\Gekt} \e f(u_\e\tx,r_{\e,k}(t)) d \sigma_x.
\end{align}
Combining the diffusion equation with the boundary condition \eqref{eq:BoundaryCondition2} and the evolution of the radii given by \eqref{eq:Derivation:Partial:r} yields the following strong formulation:
\begin{align}\label{eq:StrongForm1}
&\partial_t u_\e(t,x) - \div(D \nabla u_\e(t,x)) = f^\p(t,x)  && \textrm { in } \Oe(t),
\\\label{eq:StrongForm2}
&-D\nabla u_\e\tx u\tx  + \e \partial_t r_{\e,k}(t)\cdot n\tx  = \e f(u_\e\tx, r_{\e,k}(t))  && \textrm { in } \Gamma_{\e,k}(t),
\\\label{eq:StrongForm3}
&-D\nabla u_\e(t,x) \cdot n(t,x) = 0  && \textrm { in } \partial \Oe(t) \cap \partial \Omega,
\\\label{eq:StrongForm4}
&\partial_t r_{\e,k}(t) =  \tfrac{\e^{-N}}{c_s S_{N-1}(r_{\e,k}(t))} \int\limits_{\Gekt} f(u_\e(t,x), r_{\e,k}(t)) d \sigma_x && \textrm { for } k \in I_\e,
\end{align}
for $\Oe(t)$ given by \eqref{eq:def:Oet} and initial conditions $r_{\e}(0) = r_\e^{(0)} \in [r_\mathrm{min},r_\mathrm{max}]^{|I_\e|}$,  $ u_\e(0,\cdot_x) = u_\e^{(0)} \in L^2(\Oe(0))$.

We assume that $f^\p$ is Lipschitz continuous in every $\e$-scaled cell $\e (k + Y)$ for every $k \in I_{\e_n}$ and every $n \in \N$. Note that this does not necessarily imply $f^\p \in C(\Omega)$. 
We assume that there exists $r^{(0)}$ in $L^2(\Omega)$ such that $r^{(0)}_{\e,k(\cdot_x)} \to r^{(0)}$ in $L^2(\Omega)$, where $k_{\e}(x) \coloneqq \left(\lfloor \tfrac{x_i}{\e}\rfloor \right)_{i=1}^N$ is the index of the cell in which $x$ is located. Moreover, we assume that there exists $u_0 \in L^2(\Omega)$ such that the extension of $u_\e^{(0)}$ by $0$ to $\Omega$ two-scale converges with respect to the $L^2$-norm to $\chi_{Y_{r^{(0)}(\cdot_x)}}(\cdot_y) u_0(\cdot_x)$ for $Y_r^* \coloneqq Y \setminus \overline{B_r(x_M)}$.

\subsection{Weak formulation}
We multiply \eqref{eq:StrongForm1} by $\varphi$ and integrate over $\Oe(t)$ and $S$. Then, we
integrate the divergence term by parts and apply \eqref{eq:StrongForm2}. Thus we obtain the boundary integral $\intS\intG v_{\Gamma_{\e,k}}\tx u_\e\tx - \e f(u_\e\tx,r_{\e,k}\tx)d\sigma_x \dt$. The integration by parts of $\partial_t u \varphi$ with respect to $t$ cancels $\intS\intG v_{\Gamma_{\e,k}}\tx u_\e\tx d\sigma_x \dt$ due to the time-dependent domain $\Oe(t)$ (cf.~Reynold's transport theorem). Thus, we get \eqref{eq:WeakUnTrans1}. Furthermore, we multiply \eqref{eq:StrongForm4} by $\phi$ and integrate over $S$ which gives \eqref{eq:WeakUnTrans2}.
Altogether, we obtain the following weak form of \eqref{eq:def:Oet}, \eqref{eq:StrongForm1}--\eqref{eq:StrongForm4}:
Find $(u_\e, r_\e) \in L^2(S;H^1(\Oe(t))) \times W^{1,\infty}(S)^{|I_\e|}$ such that
\begin{align}\notag
&-\intS \int\limits_{\Oe(t)} u_\e\tx \partial_t \varphi\tx dx dt - \int\limits_{\Oe(0)}u_{\e,0}(x) \varphi(0,x) dx dt 
\\\notag
&+ \intS \int\limits_{\Oet} D \nabla u_\e\tx \cdot \nabla \varphi\tx dx dt = \intS \int\limits_{\Oet} f^\p (t,x) \varphi\tx dx dt  \\\label{eq:WeakUnTrans1}
&- \sum\limits_{k \in I_\e}\intS\int\limits_{\Gek(t)} f (u_\e\tx, r_{\e,k}(t))  \varphi\tx d \sigma_x dt
\\\label{eq:WeakUnTrans2}
&\intS \partial_t r_{\e,k}(t) \phi(t) dt = \intS \tfrac{\e^{-N} }{c_sS_{N-1}(r_{\e,k}(t)) } \int\limits_{\Gekt} \e f(u_\e\tx,r_{\e,k}(t)) d \sigma_x \phi(t) dt
\\\label{eq:WeakUnTrans3}
&r_\e(0) = r_{\e}^{(0)}
\end{align}
and \eqref{eq:def:Oet} hold 
for all $\varphi \in C^1(\overline{\bigcup\limits_{t\in S} \{t\} \times \Oe(t)})$ with $\varphi(T; \cdot_x) = 0$, all $k \in I_\e$, all $\phi \in L^1(S)^{|I_\e|}$ and all $t \in S$.
Note that $r_\e \in W^{1,\infty}(S)^{|I_\e|} \subset C^{0,1}(S)^{|I_\e|}$, which allows us to evaluate $r_{\e,k}$ pointwise in time and ensures that $\Oet$ is well defined for every $t \in S$.

\subsection{Transformation of the domain}
We transform \eqref{eq:WeakUnTrans1}--\eqref{eq:WeakUnTrans3} from $\Oe(t)$, which is given by \eqref{eq:def:Oet}, on the in time cylindrical and in space periodic domain $S \times \Oe$ with $\Oe \coloneqq \Omega \setminus \bigcup\limits_{k \in I_\e}\e\overline{ B_{r_0}(k + x_M)}$ for fixed $r_0$ with $r_\mathrm{min} \leq r_0 \leq r_\mathrm{max} $.
Thus, we can show the existence and uniqueness of a solution of \eqref{eq:WeakUnTrans1}--\eqref{eq:WeakUnTrans3} and pass to the limit $\e \to0$.
We define $\Gek \coloneqq  \partial\e B_{r_0}(k +x_M)$ for $k \in I_\e$ and $\Ge \coloneqq \bigcup_{k \in I_\e} \Gek$.

Although the geometry of $\Oet$ is already completely defined by its boundary, we need a transformation of the whole space and not only of the boundary by means of the radii, in order to apply the two-scale-transformation method.
Since $r_{\e,k} \leq r_\mathrm{max}$, the solid obstacles remain inside their respective cells so that the transformation can be defined for each $\e$-scaled cell separately using a generic transformation defined on the reference cell.

\subsubsection{Generic transformation of the reference cell}
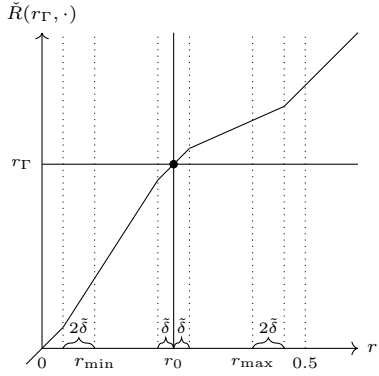
\begin{figure}[h]
\centering
\newcommand{\rzero}{0.25}
\newcommand{\rdel}{0.03}
\newcommand{\rmin}{0.1}
\newcommand{\rmax}{0.4}
\newcommand{\rGamma}{0.35}
\newcommand{\ry}{0.6}
\begin{tikzpicture}[scale = 7]
\draw[->] (0, 0) -- (0.6, 0) node[right] {\scriptsize$r$};
\draw[->] (0, 0)  node[below] {\tiny $0$}-- (0, 0.6) node[above] {\scriptsize$\check{R}(r_\Gamma, \cdot  )$};

\draw (-\rdel, -\rdel) -- (\rmin -2*\rdel , \rmin -2*\rdel);
\draw (\rmin -2*\rdel , \rmin -2*\rdel) -- (\rzero -\rdel, \rGamma-\rdel);
\draw (\rzero -\rdel, \rGamma-\rdel) -- (\rzero +\rdel, \rGamma+\rdel);
\draw (\rzero +\rdel, \rGamma+\rdel) -- (\rmax +2*\rdel , \rmax +2*\rdel);
\draw (\rmax +2*\rdel , \rmax +2*\rdel) -- (0.6, 0.6);
\draw[fill = black] (\rzero, \rGamma) circle (0.2pt);
\draw[dotted] (\rmin - 2*\rdel, 0) -- (\rmin -2*\rdel, \ry);
\draw[dotted] (\rmin, 0)  node[below] {\scriptsize $r_{\min}$} -- (\rmin , \ry);
\draw[dotted] (\rzero -\rdel, 0) -- (\rzero -\rdel, \ry);
\draw (\rzero, 0)  node[below] {\tiny $r_0$} -- (\rzero, \ry);
\draw[dotted] (\rzero +\rdel, 0) -- (\rzero +\rdel, \ry);
\draw[dotted] (\rmax, 0)  node[below] {\tiny $r_{\max}$} -- (\rmax , \ry);
\draw[dotted] (\rmax +2*\rdel, 0) -- (\rmax +2*\rdel, \ry);
\draw[dotted] (0.5, 0) node[below] {\tiny $0.5$} -- (0.5, \ry);
\draw (0.6, \rGamma) -- (0, \rGamma) node[left] {\tiny $r_\Gamma$};
\draw [decorate,decoration={brace,amplitude=4pt}](\rmin -2*\rdel, 0) -- (\rmin,0)  node [black,midway,yshift=0.3cm] {\tiny $2\tilde{\delta}$};
\draw [decorate,decoration={brace,amplitude=4pt}]
(\rzero-\rdel, 0) -- (\rzero,0) node [black,midway,yshift=0.3cm] {\tiny $\tilde{\delta}$};
\draw [decorate,decoration={brace,amplitude=4pt}]
(\rzero, 0) -- (\rzero+\rdel,0) node [black,midway,yshift=0.3cm] {\tiny $\tilde{\delta}$};
\draw [decorate,decoration={brace,amplitude=4pt}](\rmax, 0) -- (\rmax +2*\rdel,0)  node [black,midway,yshift=0.3cm] {\tiny $2\tilde{\delta}$};
\end{tikzpicture}
\caption{\label{fig:R_psi}Construction of $\check{R}$}
\end{figure}

We define the pore space of the reference cell by $\Yp \coloneqq Y_{r_0}^*$ and the interface of the reference cell by $\Gamma \coloneqq \partial B_{r_0}(x_M)$. We construct a generic cell transformation $\psi \in C^2([r_{\min}, r_{\max}] \times \overline{Y})^N$, such that
\begin{align}
\label{eq:GenericPsi_cond1}
&\psi(r_\Gamma,\Yp) = Y^*_{r_\Gamma}  \hspace{0.8cm}\textrm{ for  } r_\Gamma \in [r_{\min}, r_{\max}],
\\\label{eq:GenericPsi_cond2}
&\psi(r_\Gamma,y) = y\hspace{1.4cm}\textrm{ for } (r_\Gamma,y)  \in [r_{\min}, r_{\max}] \times \overline{Y^\p_{r_{\max}  +\delta}} \cup B_{r_{\min}  -\delta}(x_M),
\\
&\norm{\psi}{C^2([r_{\min}, r_{\max}] \times \overline{Y})} \leq C,
\\\label{eq:GenericPsi_cond_bij}
&y \mapsto \psi(r_\Gamma, y) \textrm{ is bijective from }\overline{Y} \textrm{ onto } \overline{Y},
\\\label{eq:GenericPsi_cond_det}
& \det(D_y\psi(r_\Gamma, y)) \geq c_j >0 \hspace{0.8cm}\textrm{ for  } (r_\Gamma,y) \in [r_{\min}, r_{\max}] \times \overline{Y}
\end{align}
for $\delta$ small enough. Note that due to \eqref{eq:GenericPsi_cond2}, we can glue such cell transformations $\psi(r_\Gamma,\cdot)$ for different values of $r_\Gamma$ next to each other.
\begin{figure}[h]
\centering
\begin{tikzpicture}[scale = 3]
\draw[fill=lightgray] (0,0) rectangle (1,1);
\draw[fill = white] (0.5,0.5) circle (0.45);
\CellWithBoundary{0}{0}{0.2}
\draw[fill = lightgray] (0.5,0.5) circle (0.05);

\draw[<->] (0.5,0.5) -- (0.3, 0.5);
\node at (0.4,0.58) (a) {\small$r_0$};

\node at (1.25,0.1) (b) {\small$\Yp$};
\draw[->] (b) -- (0.75,0.25);

\node at (1.50,0.65) (e1) {\small $\overline{B_{r_{\min}  -\delta}(x_M)}$};
\draw[->] (e1) -- (0.51,0.52);

\node at (1.50,0.85) (e2) {\small $\overline{Y^\p_{r_{\max}  +\delta}}$};
\draw[->] (e2) -- (0.9,0.85);
\draw[->] (1.1,0.3) -- (1.9,0.3);
\node at (1.5,0.4) (p) {\small $\psi(r, \cdot)$};
\draw[fill=lightgray] (2,0) rectangle (3,1);
\draw[fill = white] (2.5,0.5) circle (0.45);
\CellWithBoundary{2}{0}{0.35}
\draw[fill = lightgray] (2.5,0.5) circle (0.05);

\draw[<->] (2.5,0.5) -- (2.15,0.5);
\node at (2.3,0.58) (c) {\small $r$};

\node at (3.60,0.65) (f1) {\small $\psi(r,\overline{B_{r_{\min}  -\delta}(x_M)})$};
\draw[->] (f1) -- (2.51,0.52);

\node at (3.50,0.85) (f2) {\small $\psi(r,\overline{Y^\p_{r_{\max}  +\delta}})$};
\draw[->] (f2) -- (2.9,0.85);

\node at (3.6,0.25) (d) {\small $Y^*_r = \psi(r,\Yp))$};
\draw[->] (d) -- (2.8,0.08);
\end{tikzpicture}
\caption{\label{fig:psi}Generic cell transformation $\psi_(r,\cdot)$}
\end{figure}
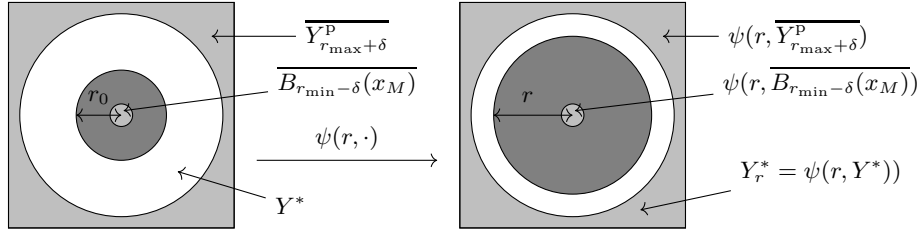
Such a generic cell transformation $\psi$ can be easily constructed
using the radial symmetry of the geometry in the reference cell.
We define 
\begin{align}\label{eq:def:psi}
\psi(r_\Gamma,y) \coloneqq  x_M +  R(r_\Gamma,\norm{y - x_M}{}) \tfrac{y - x_M}{\norm{y - x_M}{}}
\end{align}
for a smooth function $R \in C^\infty([r_{\min}, r_{\max}] \times [0,\infty))$, which scales the distance of $y$ to $x_M$ and fulfils
\newpage
\begin{align}\label{eq:cond:R1}
&R(r_\Gamma, r_0) = r_\Gamma \hspace{0.8cm}\textrm{ for } r_\Gamma \in [r_\mathrm{min},r_\mathrm{max}],
\\\label{eq:cond:R2}
&R(r_\Gamma, r) = r \hspace{1.1cm} \textrm{ for } (r_\Gamma,r) \in [r_\mathrm{min},r_\mathrm{max}] \times (\R \setminus[r_\mathrm{min} -\delta, r_\mathrm{max}+\delta]),
\\\label{eq:cond:R3}
&\norm{DR}{C^2([r_{\min}, r_{\max}]\times[0,\infty))} \leq C,
\\\label{eq:cond:R4}
&\partial_r R(r_\Gamma,r) \geq c >0  \textrm{ for }  (r_\Gamma,r) \in [r_\mathrm{min},r_\mathrm{max}] \times [0,\infty).
\end{align}
Such a mapping $R$ can be obtained by linear interpolation and smoothing (cf.~\ref{fig:R_psi}). First we define
\begin{align}
\check{R}(r_\Gamma,r) 
\coloneqq
\begin{cases}
r & \textrm{ for } r \leq r_\mathrm{min}- 2\tilde{\delta}
\\
c_1(r_\Gamma) (r-(r_\mathrm{min}- 2\tilde{\delta}) + r_\mathrm{min}- 2\tilde{\delta}& \textrm{ for }r_\mathrm{min}-2\tilde{\delta}\leq r \leq r_0-\tilde{\delta}
\\
(r-r_0) +r_\Gamma  & \textrm{ for }r_0 -\tilde{\delta} \leq r \leq r_0 + \tilde{\delta}
\\
c_2(r_\Gamma) (r-(r_\mathrm{max} + 2\tilde{\delta})) + r_\mathrm{max} +2\tilde{\delta} & \textrm{ for }r_0 +\tilde{\delta} \leq r \leq r_\mathrm{max} + 2\tilde{\delta}
\\
r & \textrm{ for } r \geq r_\mathrm{max} - 2\tilde{\delta}
\end{cases}
\end{align}
for $r_\Gamma \in [r_\mathrm{min}, r_\mathrm{max} ]$
with $c_1(r_\Gamma) \coloneqq \tfrac{r_\Gamma  - r_\mathrm{min} + \tilde{\delta}}{r_0 -r_\mathrm{min} +\tilde{\delta}}$ and
$c_2(r_\Gamma) \coloneqq \tfrac{r_\mathrm{max} - r_\Gamma +\tilde{\delta}}{r_\mathrm{max}- r_0  +\tilde{\delta}}$
and $\tilde{\delta} = \delta /3$.

Then, we define 
\begin{align}
R(r_\Gamma,r) \coloneqq \int\limits_{\R} \check{R}(r_\Gamma,s) \eta\Big(\tfrac{r-s}{\tilde{\delta}}\Big)  ds
\end{align}
for $\eta(x) \coloneqq \Big(\int\limits_\R \exp\big(\tfrac{-1}{1-|y|^2}\big) dy \Big)^{-1}  \exp\big(\tfrac{-1}{1-||^2}\big)$.
It can be shown easily that $R$ fulfils \eqref{eq:cond:R1}--\eqref{eq:cond:R4}.

We define the corresponding displacement field by $\check{\psi}(r_\Gamma, y) = \psi(r_\Gamma,y) -y$.
\subsubsection{$\e$-scaling of the transformation}

Scaling of $\psi$ by $\e$ and combining with the radii $r_{\e,k}$ for each cell gives a transformation for the $\e$-scaled porous medium:
\begin{align}\label{eq:def:psi-eps}
\psi_\e(t,x) \coloneqq [x]_{\e,Y} + \e \psi(r_{\e,k_{\e}(x)}(t), \{x\}_{\e,Y})
\end{align}
where $[x]_{\e, Y} \coloneqq \e \sum\limits_{i=1}^N k_{\e}(x)_i e_i$ is the position of the cell and $\{x\}_{\e} \coloneqq \tfrac{1}{\e} (x -[x]_{\e,Y}(x))$ is the position inside the upscaled cell.

For the corresponding displacement field, we get
\begin{align}\notag
\check{\psi}_\e\tx &\coloneqq \psi_\e\tx -x = [x]_{\e,Y} + \e \psi(r_{\e,k_{\e}(x)}(t), \{x\}_{\e,Y}) - x
\\\notag
&= [x]_{\e,Y} + \e \check{\psi}(r_{\e,k_{\e}(x)}(t), \{x\}_{\e,Y}) + \e \{x\}_{\e,Y} - x 
= \e \check{\psi}(r_{\e,k_{\e}(x)}(t), \{x\}_{\e,Y})
\end{align}
We denote the Jacobian matrix of $\psi_\e$ and its determinant by
\begin{align}\label{eq:def:Psi-eps}
\Psi_\e\tx \coloneqq D_x \psi\tx,  \ \ J_\e\tx = \det(\Psi_\e\tx).
\end{align}

Moreover, we obtain the following uniform estimates for $\psi_\e$:
\begin{lemma}[Uniform boundedness of $\psi_\e$]\label{lemma:Estimates-psi-eps}
Let $p \in[1,\infty]$, $r_{\e} \in W^{1,p}(S)^{|I_\e|}$ with $r_{\e}(t) \in [r_\mathrm{min}, r_\mathrm{max}]^{|I_\e|}$ for a.e.~$t\in S$ and let
$\psi_\e$ be defined by \eqref{eq:def:psi-eps}. Then, $\psi_\e \in W^{1,p}(S;C^1(\Oe)^N)$ and there exist constants $C,c_J,\alpha >0$ independent of $\e$ such that
\begin{align}\label{eq:Uniform_Est_psi_eps_1}
&\e^{-1}\norm{\psi_\e - \operatorname{id}_\Oe}{L^\infty(S \times \Oe)} +
\norm{\Psi_\e}{L^\infty(S \times \Oe)} +
\norm{J_\e}{L^\infty(S \times \Oe)} \leq C,
\\\label{eq:Uniform_Est_psi_eps_2}
&J_{\e}(t,x) \geq c_J,
\\\label{eq:Uniform_Est_psi_eps_3}
&\e^{-1}\norm{\partial_t \psi_\e}{L^p(S;C(\Oe))} +
\norm{\partial_t J_\e}{L^p(S;C(\Oe))} \leq C,
\\\label{eq:Uniform_Est_psi_eps_4}
&
\norm{\Psi_\e^{-1}}{L^\infty(S \times \Oe)} \leq C,
\\\label{eq:ellipticity}
&\xi^\top J_\e\tx \Psi_\e^{-1}\tx \Psi_\e^{-\top}\tx \xi\geq \alpha \norm{\xi}{}^2 
\end{align}
for a.e.~$\tx \in S\times \Oe$ and every $\xi \in \R^N$. 
\end{lemma}
\begin{proof}
The estimates \eqref{eq:Uniform_Est_psi_eps_1}--\eqref{eq:Uniform_Est_psi_eps_2} are a direct consequence of \eqref{eq:GenericPsi_cond1}--\eqref{eq:GenericPsi_cond_det} and the cell-wise construction of $\psi_\e$. The estimate \eqref{eq:Uniform_Est_psi_eps_4}--\eqref{eq:ellipticity} follow from \eqref{eq:Uniform_Est_psi_eps_1}--\eqref{eq:Uniform_Est_psi_eps_2} by simple computations.
The estimate \eqref{eq:Uniform_Est_psi_eps_3} follows with \eqref{eq:GenericPsi_cond1}--\eqref{eq:GenericPsi_cond_det} and the chain rule.
\end{proof}

Furthermore, we obtain the following uniform Lipschitz estimates for $\psi_\e$ with respect to the radii $r_\e$.
\begin{lemma}[Lipschitz regularity of $\psi_\e$]\label{lemma:Lipschitz-psi-eps}
Let $p \in [1,\infty]$ and $r_{\e,i} \in W^{1,p}(S)^{|I_\e|}$ with $r_{\e,i}(t) \in [r_\mathrm{min} , r_\mathrm{max}]^{|I_\e|}$ for a.e.~$t\in S$ and $i \in \{1,2\}$. Let $\psi_{\e,i}$ be defined by \eqref{eq:def:psi-eps} with $r_\e = r_{\e,i}$ for $i \in \{1,2\}$. Then, there exists a constant $C$ independent of $\e$ such that

\begin{align}
&\e^{-1}\norm{\psi_{\e,2}-\psi_{\e,1}}{L^\infty(S \times \Oe)} \leq C\norm{r_{\e,2}- r_{\e,1}}{L^\infty(S)},
\\
&\norm{\Psi_{\e,2}-\Psi_{\e,1}}{L^\infty(S \times \Oe)} + \norm{J_{\e,2}-J_{\e,1}}{L^\infty(S \times \Oe)}\leq C\norm{r_{\e,2}- r_{\e,1}}{L^\infty(S)},
\\
&\norm{\Psi_{\e,2}^{-1}-\Psi_{\e,1}^{-1}}{L^\infty(S \times \Oe)} \leq C\norm{r_{\e,2}- r_{\e,1}}{L^\infty(S)}, \\
&\e^{-1}\norm{\partial_t (\psi_{\e,2} - \psi_{\e,1})}{L^p(S\times \Oe)} +
\norm{\partial_t (J_{\e,2} - J_{\e,1})}{L^p(S\times \Oe)} \leq C \norm{\partial_t (r_{\e,2} - r_{\e,1})}{L^p(S\times \Oe)}.
\end{align}
\end{lemma}
\begin{proof}
Lemma \ref{lemma:Lipschitz-psi-eps} can be proven by similar computations as in the proof of Lemma \ref{lemma:Estimates-psi-eps}.
\end{proof}

\subsection{Transformation of the weak form}
Using the diffeomorphism $\psi_\e$, which is defined in \eqref{eq:def:psi-eps}, we define
$\hat{f}^\p_\e(t,x)\coloneqq f^\p(t,\psi_\e\tx)$ and note that
Lemma \ref{lemma:Estimates-psi-eps} implies the uniform estimate for $\hat{f}^\p_\e$ by
\begin{align}
\norm{\hat{f}_\e^\p}{S \times \Oe}^2 &= \int\limits_{S\times \Oe} f^\p(t,\psi_\e\tx)^2 \dxt = \intS \int\limits_\Oet J_\e^{-1}(t,\psi_\e^{-1}\tx) f_\e^\p\tx^2 \dxt
\\
&\leq c_J^{-1} \intS \int\limits_\Oet f_\e^\p\tx^2 \dxt \leq C \norm{f^\p}{S \times \Omega}^2 
\end{align}
We define $A_\e \coloneqq J_\e \Psi_\e^{-1} D \Psi_\e^{-\top}$ and $B_\e \coloneqq J_\e \Psi_\e^{-1} \partial_t \psi_\e$.
Then, we transform the weak form \eqref{eq:def:Oet},\eqref{eq:WeakUnTrans1}--\eqref{eq:WeakUnTrans3} into the following equivalent weak form:

Find $(\hue, r_\e) \in L^2(S;H^1(\Oe)) \times W^{1,\infty}(S)^{|I_\e|}$ such that $\partial_t (J_\e u_\e)  \in L^2(S;H^1(\Oe)')$ and
\begin{align}\notag
&\langle \partial_t(J_\e(t) \hue(t),\varphi \rangle_{\Oe}
+ (A_\e(t) \nabla \hue(t), \nabla \varphi(t))_{\Oe}
+ (B_\e(t) \hue(t), \nabla \varphi) 
\\\label{eq:WeakTrans1}
&= (J_\e(t)\hat{f}^\p_\e(t), \varphi)_\Oe 
-\sum\limits_{k \in I_\e}\tfrac{r_{\e,k}^{n-1}(t)}{ r_0^{n-1}} (   \e f_\e (\hue(t), r_{\e,k}(t))  \varphi(t))_\Gek
\\\label{eq:WeakTrans2}
&\intS \partial_t r_{\e,k}(t) \phi(t) dt = \intS \tfrac{\e^{-N} }{c_sS_{N-1}(r_0) } \int\limits_{\Gek} \e f(\hue\tx,r_{\e,k}(t)) d \sigma_x \phi(t) dt,
\\\label{eq:WeakTrans3}
&r_\e(0) = r_{\e}^{(0)}, \ \ \hue(0) = \hue^{(0)} \coloneqq u_\e^{(0)}\circ \psi_0^{-1}(0)
\end{align}
for a.e.~$t \in S$ all $\varphi \in H^1(\Oe)$, all $k \in I_\e$ and all $\phi \in L^1(S)^{|I_\e|}$, where $\psi_\e$ depends on $r_\e$ and is defined by \eqref{eq:def:psi-eps} and $\Psi_\e, J_\e $ by \eqref{eq:def:Psi-eps}.

\begin{lemma}
Let $\psi_\e, \Psi_\e, J_\e$ be given by \eqref{eq:def:psi-eps} and \eqref{eq:def:Psi-eps}, respectively. Then, 
$(u_\e,r_\e)$ is a solution of \eqref{eq:def:Oet}, \eqref{eq:WeakUnTrans1}--\eqref{eq:WeakUnTrans3} if and only if $\hue = u_\e(\cdot_t,\psi_\e(\cdot_t, \cdot_x))$ is a solution of \eqref{eq:def:psi-eps}--\eqref{eq:def:Psi-eps}, \eqref{eq:WeakTrans1}--\eqref{eq:WeakTrans3}.
\end{lemma}
\begin{proof}
The proof follows by a simple transformation and the density of $C^1(S\times\Oe) \subset L^2(S;H^1(\Oe))$.
\end{proof}

\section{Existence and uniform a priori estimates}
For the existence proof, we combine a fixed-point argumentation with the theory of monotone operators from \cite{Sho97}.
\begin{defi}[Monotone operator]
A function $\A : V \rightarrow V'$ is \textit{monotone} if $(\A(u) -\A(v) , u-v)_V \geq 0$ for every $u,v \in V$.
\end{defi}
\begin{defi}[Family of regular operators]
Let $W$ be a separable Hilbert space. A family of operators $\{B(t) \mid t \in \overline{S}$ with $B(t) \in L(W,W')$ for each $t \in \overline{S}$ and $B(\cdot)(u)(v) \in L^\infty(S)$ for each pair $u,v \in W$ is called \textit{regular} if for each pair $u,v \in W$, the function $B(\cdot)u (v)$ is absolutely continuous on $\overline{S}$ and there is a $K \in L^1(S)$ such that 
\begin{align}
|\tfrac{d}{dt} B(t) u(v) | \leq K(t) \norm{u}{W} \norm{u}{W}
\end{align}
for every $u,v \in W$ and for a.e.~$t\in \overline{S}$.
\end{defi}
The monotone operator theory gives the following existence result for degenerate parabolic equations (cf.~\cite{Sho97}).
\begin{thm}\label{thm:existenceDegenerateParabolic}
Let $V$ be a separable Hilbert space. Suppose that $W$ is a Hilbert space containing $V$ with dense and continuous injection $V  \xhookrightarrow{} W$. Let $\V \coloneqq L^2(S;V)$ and $\W \coloneqq L^2(S;W)$. We assume that for every $t \in \overline{S}$ there are given operators $\A(t) \in L(V,V')$ and $\B(t) \in L(W,W')$ such that $\A(\cdot) u(v) \in L^\infty(S)$ for each pair $u,v \in V$ and $\B(\cdot) u(v) \in L^\infty(S)$ for each pair $u,v \in W$. 

In addition, we assume that $\{\B(t) \mid t \in \overline{S}\}$ is a regular family of self-adjoint operators, $\B(0)$ is monotone and there are numbers $\lambda, c >0$ such that 
\begin{align}\label{eq:existenceDegenerateParabolicCoercivity}
2 \A(t)v(v) + \lambda B(t)v(v) + \B'(t)v(v) \geq c \norm{v}{V} &&\textrm{for all } V \in V \textrm{ and all } t \in \overline{S}.
\end{align}
Then, for given $u^{(0)} \in W$ and $f \in L^2(0,T;V')$ there exists $u \in \V$ such that
\begin{align}\label{eq:degenerateAbstractEq}
\tfrac{d}{dt}(\B(t)u(t)) + \A(t)u(t) = f(t)   \textrm{ in } \V', \textrm{ with } (\B u)(0) = \B(0)u_0.
\end{align}
\end{thm}
Combining Theorem \ref{thm:existenceDegenerateParabolic} with a fixed point argumentation allows us to derive the existence and uniqueness of the solution of the system \eqref{eq:def:psi-eps}--\eqref{eq:def:Psi-eps}, \eqref{eq:WeakTrans1}--\eqref{eq:WeakTrans3}.
\begin{thm}\label{thm:Existence-ue}
There exists a unique solution $(\hat{u}_\e,r_\e) \in L^2(S;H^1(\Omega))\times W^{1,\infty}(S)^{|I_\e|}$ with $\partial_t(J_\e \hat{u}_\e), \partial_t\hat{u}_\e \in L^2(S;H^1(\Oe)')$ of the system \eqref{eq:def:psi-eps}--\eqref{eq:def:Psi-eps}, \eqref{eq:WeakTrans1}--\eqref{eq:WeakTrans3} and thus $\hue \in C^0(\overline{S};L^2(\Oe))$. Moreover, the following uniform estimates hold
\begin{align}\label{eq:estimate-u-eps}
&\norm{\hue}{C^0(\overline{S};L^2(\Oe))} + \norm{\nabla \hue}{L^2(S \times \Oe)} \leq C,
\\\label{eq:estimate-r-eps}
&r_{\e,k}(t) \in [r_\mathrm{min} , r_\mathrm{max}]^{|I_\e|} \textrm{ for every } t \in \overline{S},
\\\label{eq:estimate-dt_r-eps}
&\norm{\partial_t r_{\e,k}}{L^\infty(S)} \leq C_f c_s^{-1}  \textrm{ for every } k \in I_\e.
\end{align} 
\end{thm}
\begin{proof}
In order to show the existence and uniqueness of the solution, we divide $S$ in finitely many subintervals $S_i \coloneqq (t_i,t_{i+1})$ with $0 = t_0 <t_1 < \dots < t_n = T$ for $i \in \{0, \dots, N_\e\}$ and $N_\e$ large enough.
Then, we show iteratively that there exists a unique solution $(\hue|_{S_i}, r_\e|_{S_i} )\in L^2(S_i;H^1(\Oe)) \times W^{1,\infty}(S_i)^{|I_\e|}$ with $\partial_t (J_\e \hue), \partial_t \hue \in L^2(S_i;H^1(\Oe)')$ such that \begin{align}\notag
\\\notag
&\int\limits_{t_i}^{t_{i+1}} \langle \partial_t(J_\e(\tau) \hue|_{S_i}(\tau),\varphi(\tau) \rangle_{\Oe} \dd \tau
+ (A_\e \nabla \hue|_{S_i}, \nabla \varphi)_{(t_i,t_{i+1}) \times \Oe}
+ (B_\e \hue|_{S_i}, \nabla \varphi)_{(t_i,t_{i+1}) \times \Oe} 
\\\label{eq:weakFormSi1}
&= (J_\e\hat{f}^\p_\e, \varphi)_{(t_i,t_{i+1}) \times \Oe}  
-\sum\limits_{k \in I_\e} \left( \tfrac{r_{\e,k}^{n-1}}{ r_0^{n-1}}   \e f_\e (\hue|_{S_i}, r_{\e,k}),  \varphi\right)_{(t_i,t_{i+1}) \times \Gek}
\\\label{eq:weakFormSi2}
&\int\limits_{t_i}^{t_{i+1}} \partial_t r_{\e,k}(t) \phi(t) dt = \int\limits_{t_i}^{t_{i+1}}\tfrac{\e^{-N} }{c_sS_{N-1}(r_0) } \int\limits_{\Gek} \e f(\hue|_{S_i}\tx,r_{\e,k}(t)) d \sigma_x \phi(t) dt,
\end{align}
holds for every $(\varphi,\phi) \in L^2(S_i;H^1(\Oe))\times L^2(S_i)^{|I_\e|}$ and the initial condition $(\hue|_{S_i}(t_i), r_\e|_{S_i}(t_i) ) = (\hat{u}_\e^{(t_i)}, r_\e^{(t_i)})$ is fulfilled. For $i \geq 1$, the initial values  are defined by means of the solution on the previous time interval, i.e.~$ (\hat{u}_\e^{(t_i)}, r_\e^{(t_i)}) \coloneqq (\hue|_{S_{i-1}}(t_i), r_\e|_{S_{i-1}}(t_i))$. Then, we get the solution $(\hue,r_\e)$ for the whole interval $S$ by concatenating the solutions.

First, we choose $t_1$ small enough such that we can apply Lemma \ref{lemma:Existence_On_Si}. Then, we get a solution $(\hue|_{S_0}, r_\e|_{S_0}) \in L^2(S_0;H^1(\Omega)) \times W^{1,\infty}(S_0)^{|I_\e}$ with $\partial_t \hue|_{S_0} \in L^2(S_i;H^1(\Oe)')$. 
Now, we proceed inductively. We assume that we have a unique solution $(\hue|_{(0, t_i)}, r_{\e}|_{(0,t_i)})$ of \eqref{eq:def:psi-eps}--\eqref{eq:def:Psi-eps}, \eqref{eq:WeakTrans1}--\eqref{eq:WeakTrans3} for the time interval $(0,t_i)$ instead of $S$. Then, we claim that there exists also an unique solution on the time interval $(0,t_{i+1})$ where $t_{i+1} -t_i \geq  \sigma_\e >0$ for a constant $\sigma_\e$ which depends nether on the iteration number $i$ nor on the exact time $t_i$ as long as $t_i \leq T$. Hence, we obtain after finitely many steps a solution for the whole interval $S$.
In order to show this uniform bound $\sigma_\e$, we use Lemma \ref{lemma:Existence_On_Si} and note that we have only to show that
\begin{align}\label{eq:UniformBound_r_eps_init}
&r_\e|_{(0,t_i)}(t_i) \in [r_{\min},r_{\max}]^{|I_\e|},
\\\label{eq:UniformBound_u_eps_init}
&\norm{\hue|_{(0,t_i)}(t_i)}{\Oe}\leq K,
\end{align}
for a constant $K$ which is independent on the iteration number $i$ and the time $t_i\leq T$. Then, we can construct the solution on $(t_i, t_{i+1})$ with Lemma \ref{lemma:Existence_On_Si} and can concatenate it with the solution on $(0,t_{i})$.
The estimate \eqref{eq:UniformBound_r_eps_init} follows directly from Lemma \ref{lemma:Existence_On_Si} since $r_\e|_{(0,t_i)}$ was constructed by Lemma \ref{lemma:Existence_On_Si}.
The estimate \eqref{eq:UniformBound_u_eps_init} can be derived like the estimates \eqref{eq:weakFormSi1_Tested_ue}--\eqref{eq:est:hue_Si} but applied on $\hue|_{(0,t_i)}$.
The crucial point is that the constant in \eqref{eq:est:hue_Si} does not depend on $t$ as long as $t\leq T$. It depends only on the initial value. Since we do not apply the estimates iteratively on the interval $(t_l,t_{l+1})$ for $l \in \{0,\dots , i-1\}$ but only once on the whole interval $(0,t_i)$, we do not have to take care if the initial values multiply in a bad manner. However, we have to note that these estimates formally only bound $\norm{\hue}{L^\infty((0,t_i);L^2(\Oe))}$ uniformly. In order to get the  uniform bound not only for a.e.~$t\in(0,t_i)$ but also for $t_i$, we use the following argument. Since $\partial_t \hue \in L^2((0,t_i);H^1(\Oe)')$, the Lemma of Lions-Aubin gives $\hue\in C(\overline{(0,t_i)};L^2(\Oe))$ and since $\norm{\hue}{L^\infty((0,t_i);L^2(\Oe))} = \norm{\hue}{C(\overline{(0,t_i)};L^2(\Oe))}$, we get the uniform bound for $\norm{\hue|_{(0,t_i)}(t_i)}{\Oe}$.

Moreover, we note that the estimates \eqref{eq:weakFormSi1_Tested_ue}--\eqref{eq:est:hue_Si}  do not depend on $\e$.
In fact, an $\e$-dependency would not be a problem for the the proof of the existence and uniqueness of $\hue$ on the whole time interval. However, due to their $\e$-independency, they give us immediately the uniform bound \eqref{eq:estimate-u-eps} since the initial values $\hue^{(0)}$ are uniformly bounded.
\end{proof}

\begin{lemma}\label{lemma:Existence_On_Si}
Let $\e >0$ and $S_i = (t_i,t_{i+1})$ with $0\leq t_i < t_{i+1}\leq T$. Then, for every $K>0$, exists a constant $\sigma_{\e,K}>0$, which depends only on $\e$ and $K$, such that  \eqref{eq:def:psi-eps}--\eqref{eq:def:Psi-eps}, 
\eqref{eq:weakFormSi1}--\eqref{eq:weakFormSi2}, has a 
unique solution $(\hue, r_\e) \in L^2(S_i;H^1(\Oe)) \times W^{1,\infty}(S)^{|I_\e|}$ 
with $\partial_t (J_\e \hat{u}_\e), \partial_t \hat{u}_\e \in L^2(S;H^1(\Oe)')$, $\hue(t_i) = \hat{u}_\e^{(t_i)}$ and $r_\e(t_i) = r_{\e}^{(t_i)}$ for arbitrary 
$\hat{u}_\e^{(t_i)} \in L^2(\Oe)$ and
$r_{\e}^{(t_i)} \in [r_\mathrm{min},r_\mathrm{max}]^{|I_\e|}$,
if $\norm{\hat{u}_\e^{(t_i)}}{\Oe} \leq K$ and
$|S_i|\leq \max \{1,\sigma_{\e,K}\}$.
Moreover, $r(t) \in [r_\mathrm{min},r_\mathrm{max}]^{|I_\e|}$ and $|\partial_t r(t)| \leq C_f c_s^{-1}$ for a.e.~$t \in S_i$.
\end{lemma}

\begin{proof}
We show the existence and uniqueness by means of a fixed-point argumentation for $\hue \in L^2(S_i;H^1(\Oe))$ with the fixed-point operator $L_\e : L^2(S_i;H^1(\Oe)) \rightarrow L^2(S_i;H^1(\Oe))$.
First, $L_\e$ inserts a given function $\zeta$ into the right-hand side of \eqref{eq:weakFormSi1}--\eqref{eq:weakFormSi2}, which yields
\begin{align}\notag
&\int\limits_{t_i}^{t_{i+1}} \langle \partial_t(J_\e(t) \hue(t),\varphi(t) \rangle_{\Oe} \dt
+ (A_\e \nabla \hue, \nabla \varphi)_{(t_i,t_{i+1}) \times \Oe}
+ (B_\e \hue, \nabla \varphi)_{(t_i,t_{i+1}) \times \Oe} 
\\\label{eq:FixedPointUe}
&= (J_\e\hat{f}^\p_\e, \varphi)_{(t_i,t_{i+1}) \times \Oe}  
-\sum\limits_{k \in I_\e} \left( \tfrac{r_{\e,k}^{n-1}}{ r_0^{n-1}}   \e f(\zeta, r_{\e,k}),  \varphi\right)_{(t_i,t_{i+1}) \times \Gek},
\\\label{eq:FixedPointRe}
&\int\limits_{t_i}^{t_{i+1}} \partial_t r_{\e,k}(t) \phi(t) dt = \int\limits_{t_i}^{t_{i+1}}\tfrac{\e^{-N} }{c_sS_{N-1}(r_0) } \int\limits_{\Gek} \e f(\zeta\tx,r_{\e,k}(t)) d \sigma_x \phi(t) dt.
\end{align} 
Then it solves \eqref{eq:FixedPointRe} for $r_\e$. This $r_\e$ gives $\psi_\e, \Psi_\e, J_\e$ via \eqref{eq:def:psi-eps}--\eqref{eq:def:Psi-eps} for \eqref{eq:FixedPointUe}. Then, $L_\e(\zeta) \coloneqq\hue$ where $\hue$ is the solution of  \eqref{eq:FixedPointUe}.

In order to show that $L_\e$ is well defined and is a contraction, we rewrite $L_\e(\hue)$ by means of the following both operators.
Let $V_{r,\e}(S_i) \coloneqq\{r \in  W^{1,2}(S_i)^{|I_\e|} \mid r(t) \in [r_\mathrm{min},r_\mathrm{max}]^{|I_\e|} \textrm{ and } |\partial_t r(t)| \leq C_f c_s^{-1} \textrm{ for a.e.~}t \in S_i\}$.
We define $L_{\e,1} : L^2(S_i;H^1(\Oe)) \rightarrow V_{r,\e}(S_i)$
as the solution operator of \eqref{eq:FixedPointRe}, i.e.~$L_{\e,1}(\zeta) \coloneqq r_\e$, where $r_\e \in V_{r,\e}(S_i)$ is the solution of \eqref{eq:FixedPointUe} for every $k \in I_\e$ and every $\phi \in L^2(S_i)$ with initial condition $r_\e(t_i) = r_{\e}^{(t_i)}$. 
Moreover, we define $L_{\e,2} : L^2(S_i;H^1(\Oe)) \times V_{r,\e}(S_i) \rightarrow L^2(S_i;H^1(\Oe))$
by $L_{\e,2}(\hat{\zeta}_\e, r_\e) \coloneqq \hue$, where $\hue$ is the solution of
\eqref{eq:FixedPointUe} for every $\varphi \in L^2(S;H^1(\Oe))$ with initial condition $\hue(t_i) = \hat{u}_\e^{(t_i)}$.
Hence, we get $L_\e(\zeta) = L_{\e,2}(\zeta, L_{\e,1}(\zeta))$. 

Note, that $\hue$ is a fixed point of $L_\e$ with $\partial_t (J_\e \hue), \partial_t \hue \in L^2(S;H^1(\Oe)')$ and $r_\e = L_{\e,i}{\hue}$ with $\partial_t r_\e \in L^\infty(S_i)^{|I_\e|}$ if and only if $(\hue,r_\e)$ solves \eqref{eq:weakFormSi1}--\eqref{eq:weakFormSi2}.

Hence, it is sufficient to show, that $L_\e$ has a unique fixed point. First, we show, that $L_{\e,1}$ is well defined and Lipschitz continuous. Then, we do the same for $L_{\e,2}$. Thereby, we show that the Lipschitz constants of $L_{\e,1}$ and $L_{\e,2}$ tend to zero for $|S_i| \to 0$. Thus, we obtain that $L_\e$ is a contraction for $|S_i|$ small enough and the contraction theorem gives the existence and uniqueness of a fixed point of $L_\e$.

\paragraph{$\bullet$ $L_{\e,1}$ is well defined}

Since $(t,r) \mapsto \tfrac{\e^{-N}}{c_sS_{N-1}(r_0) } \int_{\Gek} \e f(\zeta\tx,r) d \sigma_x$ is globally Lipschitz continuous with respect to $r$ and measurable with respect to $t$ if $\zeta \in L^2(S;H^1(\Oe))$, Carath\'eodory's existence theorem yields the existence and uniqueness of a solution $r_{\e} \in W^{1,1}(S_i)$ of \eqref{eq:weakFormSi2}. Moreover,  the Assumption \eqref{eq:f>0}--\eqref{eq:f<0} ensure that $L_{\e,1}(r_{\e,i})(t) \in [r_\mathrm{min},r_\mathrm{max}]^{|I_\e|}$ for a.e.~$t\in S_i$ and \eqref{eq:f-bounded} that $\left|\tfrac{\e^{-N}}{c_sS_{N-1}(r_0) } \int_{\Gek} \e f(\zeta\tx,r) d \sigma_x\right| \leq C_fc_s^{-1}$. Thus $L_{\e,1}$ is well defined with $L_{\e,1}(\zeta) = r_\e   \in V_{r,\e}(S_i)$.

\paragraph{$\bullet$ Lipschitz estimate of $L_{\e,1}$}
For the Lipschitz estimate, let $\zeta_1, \zeta_2 \in L^2(S;H^1(\Oe))$. We define $r_{\e,i} \coloneqq L_{\e,1}(\zeta_i)$ for $i \in \{1,2\}$ and test \eqref{eq:FixedPointRe} for $\zeta = \zeta_i$ for $i \in \{1,2\}$ with $\chi_{(t_i,t)}(r_{\e,1,k} - r_{\e,2,k})$ for $t\in (t_i,t_{i+1})$. We subtract both equations. Then, we obtain with the Lipschitz condition \eqref{eq:f-Lipschitz} of $f$, the Young and the Cauchy--Schwarz inequalities
\begin{align*}
&\tfrac{1}{2}|r_{\e,1,k}(t) - r_{\e,2,k}(t)|^2 = ( \partial_t(r_{\e,1,k} - r_{\e,2,k}),r_{\e,1,k} - r_{\e,2,k} )_{(t_i,t)}
\\ 
&= \tfrac{\e^{-N} }{c_sS_{N-1}(r_0)} \e (
f(\hat{u}_{\e,1},r_{\e,1,k}) -f(\hat{u}_{\e,2},r_{\e,2,k}), r_{\e,1,k} - r_{\e,2,k} )_{(t_i ,t) \times \Gek}
\\
&\leq 
\tfrac{\e^{-N} }{c_sS_{N-1}(r_0)} \e \norm{
C_{L_f}
(\hat{u}_{\e,1}- \hat{u}_{\e,2}) + (r_{\e,1,k} -r_{\e,2,k})}{(t_i ,t) \times \Gek} 
\norm{r_{\e,1,k} - r_{\e,2,k}}{(t_i ,t) \times \Gek}
\\
&\leq
C_\e \norm{r_{\e,1,k} -r_{\e,2,k}}{(t_{i},t)}^2 + C_\e \norm{\zeta_1 - \zeta_2}{(t_i , t) \times \Gek}^2
\\
&\leq
C_\e  \norm{r_{\e,1,k} -r_{\e,2,k}}{(t_{i},t)}^2 + C_\e  \norm{\zeta_1 - \zeta_2}{S_i\times \Ge}^2
\end{align*}
After collecting all the constants and applying the Lemma of Gronwall, we get
\begin{align}\label{eq:L(r)<u}
|r_{\e,1,k}(t) - r_{\e,2,k}(t)|^2
\leq
C_\e \norm{\zeta_1 - \zeta_2}{S_i\times \Ge}^2 
\end{align}
for every $t \in S_i$.

Then, we test \eqref{eq:FixedPointRe} for $\zeta = \zeta_i$ for $i \in \{1,2\}$ with $\partial_t r_{\e,1,k} - \partial_t r_{\e,2,k}$ and use again the Lipschitz condition \eqref{eq:f-Lipschitz}:
\begin{align}\notag
&\norm{\partial_t r_{\e,1,k} - \partial_t r_{\e,2,k}}{S_i}^2
\\\notag
&= 
\tfrac{\e^{-N} }{c_sS_{N-1}(r_0)} \e (f(\zeta_1, r_{\e,1,k}) - f(\zeta_2, r_{\e,2,k}) ,  \partial_t r_{\e,1,k}- \partial_t r_{\e,2,k})_{S_i \times \Gek}
\\\notag
&\leq
C_\e\norm{
C_{L_f} (\zeta_1- \zeta_2) + (r_{\e,1,k} -r_{\e,2,k})}{S_i \times \Gek} 
\norm{\partial_t r_{\e,1,k} - \partial_t r_{\e,2,k}}{S_i}
\\\label{eq:dt-L(r)<u}
&\leq 
C_\e (\norm{\zeta_1- \zeta_2}{S_i \times \Gek}  + \norm{r_{\e,1,k} -r_{\e,2,k}}{S_i})
\norm{\partial_t r_{\e,1,k} - \partial_t r_{\e,2,k}}{S_i}.
\end{align}
Inserting \eqref{eq:L(r)<u} in \eqref{eq:dt-L(r)<u} and employing the continuity of the trace operator for $\Gek$ yields
\begin{align}
&\norm{\partial_t r_{\e,1,k} - \partial_t r_{\e,2,k}}{S_i}
\leq
C_\e \norm{\zeta_1 -\zeta_2}{S \times \Gek} 
\label{eq:dt-L(r)<dxu}
\leq
C_\e \norm{\zeta_1 -\zeta_2}{L^2(S_i;H^1(\Oe)}.
\end{align}
The fact that $r_\e \in V_{r,\e}(S_i)$ implies $\norm{\partial_t r_{\e,k}}{L^\infty(S_i)} \leq C_f c_s^{-1}$. Thus, we get with the H\"older inequality
\begin{align}\notag
&\norm{\partial_t r_{\e,1,k} - \partial_t r_{\e,2,k}}{S_i}
\leq
\norm{1}{S_i}\norm{(\partial_t r_{\e,1,k} - \partial_t r_{\e,2,k})^2}{S_i}
\\\notag
&\leq
\norm{1}{S_i} \norm{\partial_t r_{\e,1,k} - \partial_t r_{\e,2,k}}{L^\infty(S_i)} \norm{\partial_t r_{\e,1,k} - \partial_t r_{\e,2,k}}{S_i}
\\\label{eq:LipEstL1a}
&\leq 
\sqrt{|S_i|} C_\e \norm{\zeta_1 -\zeta_2}{L^2(S_i;H^1(\Oe)} 
\end{align}
Moreover, we can conclude with the fundamental theorem of calculus and the H\"older inequality fo every $t \in S_i$:
\begin{align}\notag
&|r_{\e,1,k} -r_{\e,2,k}(t)|
=
\int\limits_{t_i}^t \partial_t (r_{\e,1,k} -r_{\e,2,k})(\tau) \dd \tau
\leq \norm{1}{S_i}
\norm{\partial_t r_{\e,1,k} - \partial_t r_{\e,2,k}}{S_i}
\\\label{eq:LipEstL1b}
&\leq |S_i| C_\e \norm{\zeta_1 -\zeta_2}{L^2(S_i;H^1(\Oe)} 
\end{align}

\paragraph{$\bullet$ $L_{\e,2}$ is well defined}

First, we show the existence of a solution $\hue \in L^2(S_i;H^1(\Oe))$ with $\partial_t(J_\e \hue) \in L^2(S;H^1(\Oe)')$ of \eqref{eq:FixedPointUe} using Theorem \ref{thm:existenceDegenerateParabolic}. 
With  the regularity of $J_\e$, we can conclude $\partial_t \hue \in L^2(S;H^{1}(\Oe)')$.
Testing \eqref{eq:FixedPointUe} with $\hue$ shows the uniqueness of the solution of \eqref{eq:FixedPointUe} and thus that $\hue = L_{\e,2}(\zeta,r_\e)$ is well defined for every $\zeta \in L^2(S;H^1(\Oe))$.

Using the setting of Theorem \ref{thm:existenceDegenerateParabolic}, we set $V= H^1(\Oe)$ and $W = L^2(\Oe)$. 
Let $\psi_\e, \Psi_\e$ and $J_\e$ be given by \eqref{eq:def:psi-eps}--\eqref{eq:def:Psi-eps}.
For each $t \in[t_i, t_i+1]$ and $u,v \in V$, we define $\A_\e(t) : V \rightarrow V'$ by $
(\A_\e(t) u)(v) \coloneqq (A_\e(t) \nabla u, \nabla v)_\Oe + (B_\e(t) u, \nabla v)_\Oe
$.
For each $t \in[t_i, t_{i+1}]$ and $u,v \in W$, we define $\B_\e(t) : W \rightarrow W'$ by $
(\B_\e(t) u)(v) \coloneqq (J_\e(t)u,v)_\Oe$.
For $\zeta,v \in \V$, we define $f_\e(\zeta; \cdot) : \V \rightarrow \R$ by
\begin{align*}
f_\e(\zeta;v) \coloneqq (J_\e \hat{f}_\e^\p, v)_{S_i\times \Oe}
-\sum\limits_{k \in I_\e} ( \tfrac{r_{\e,k}^{n-1}}{r_0^{n-1}}  f(\zeta,r_{\e,k}), v )_{S_i\times \Gek}.
\end{align*}
In order to apply Theorem \ref{thm:existenceDegenerateParabolic}, we verify its assumption in the following.
The Lipschitz regularity of $f$ and the continuous embedding $H^1(\Oe)\xhookrightarrow{} L^2(\Gek)$ ensure that $f_\e(\zeta;\cdot) \in \V'$ for every $\zeta \in L^2(S;H^1(\Oe))$. Moreover, it is clear that $\A_\e(t) \in L(V,V')$, $\B_\e(t) \in L(W,W')$ for every $t\in [t_i,t_{i+1}]$. 
Since $r_\e \in V_{r,\e}(S_i)$, we can conclude with Lemma \ref{lemma:Estimates-psi-eps} that $\A_\e(\cdot)u(v) \in L^\infty(S)$ for every pair $u,v \in V$ and $\B_\e(\cdot)u(v) \in L^\infty(S)$ for every pair $u,v \in W$.
Furthermore, it is clear that $\{\B_\e(t) \mid t \in [t_i,t_{i+1}] \}$ is a family of self-adjoint operators. 
From Lemma \ref{lemma:Estimates-psi-eps}, we get the time regularity of $J_\e$ which can be transferred on $\B_\e$ so that $\{\B_\e(t) \mid t \in [t_i,t_{i+1}] \}$ is a family of regular operators.
Using the uniform boundedness of $J_\e$ from below given by Lemma \ref{lemma:Estimates-psi-eps}, we get that $\B(0)$ is monotone.
It remains to show the estimate \eqref{eq:existenceDegenerateParabolicCoercivity}.
Using the coercivity of $J_\e \Psi_\e^{-1} \Psi_\e^{-\top}$ given by Lemma \ref{lemma:Estimates-psi-eps}, we obtain for every $v \in H^1(\Oe)$ and every $t\in \overline{S}$ 
\begin{align}\label{eq:Coercivity_A_e}
(A_\e(t) \nabla vy \nabla v)_\Oe \geq \alpha \norm{\nabla v}{\Oe}.
\end{align}
Using the estimates on $\Psi_\e, J_\e$ and $\partial_t \psi_\e$ of Lemma \ref{lemma:Estimates-psi-eps} as well as the H\"older and Young inequalities, we get for every $\delta >0$ a constant $C_\delta$ such that for every $v \in H^1(\Oe)$ and every $t\in \overline{S}$ 
\begin{align}\label{eq:Estimate_B_e}
-(B_\e(t) v,  \nabla v)_\Oe \leq C \norm{v}{\Oe}\norm{\nabla v}{\Oe} \leq C_\delta \norm{v}{\Oe}^2 + \delta\norm{\nabla v}{\Oe}^2 
\end{align}
Combing \eqref{eq:Coercivity_A_e}--\eqref{eq:Estimate_B_e} with the definition of $\A_\e(t)$ yields for $\delta = \alpha/2$
\begin{align}\label{eq:Coercivity_A}
\A_\e(t)v(v) = (A_\e(t) \nabla v, \nabla v)_\Oe + (B_\e(t) v,  \nabla v)_\Oe \geq \alpha/2 \norm{\nabla v }{\Oe}^2 - C_{\alpha/2}\norm{v}{\Oe}^2
\end{align}
The estimate on $J_\e$ from below implies
\begin{align}
\B_\e(t) v(v) \geq c_J \norm{v}{\Oe}^2
\end{align}
and the boundedness of $\norm{\partial_t r_\e}{L^\infty(S_i)} \leq C$  together with Lemma \ref{lemma:Estimates-psi-eps} gives
\begin{align}
-\B'(t)v(v) = (\partial_t J_\e(t) v, v)_\Oe\leq C \norm{v}{\Oe}^2
\end{align}
Thus, we get
\begin{align}\label{eq:Coercivity_B'+B}
\lambda \B_\e(t) v(v)  + \B'(t)v(v) \geq (\lambda c_J -C) \norm{v}{\Oe}^2.
\end{align}
Combining  \eqref{eq:Coercivity_A}--\eqref{eq:Coercivity_B'+B} for $\lambda = (\alpha/2 + C - C_{\alpha/2})/c_J$ gives
\eqref{eq:existenceDegenerateParabolicCoercivity}.
Thus, we have shown all that all prerequisites of Theorem \ref{thm:existenceDegenerateParabolic} are fulfilled and we get a solution $\hue \in L^2(S_i;H^1(\Oe))$ with $\partial_t (J_\e \hue )\in L^2(S_i;H^1(\Oe)')$. Then, the regularity of $J_\e$ implies that $\partial_t \hue = \langle \partial_t(J_\e \hue), \cdot \rangle_\Oe - (\partial_t J_\e \hue , \cdot )_\Oe \in L^2(S;H^1(\Oe)')$.

In order to show that $L_{\e,2}$ is well defined, it remains to show the uniqueness of the solution of \eqref{eq:FixedPointUe}. 
Due to the linearity of the equation \eqref{eq:FixedPointUe}, it is sufficient to show that $\hue= 0$, if $\hat{u}_\e^{(t_i)} = 0$, $\hat{f}^\p_\e=0$ and $f = 0$.
Therefore, we test \eqref{eq:FixedPointUe} with the solution $\chi_{(t_i,t)}\hue$ for $t \in S_i$, which yields
\begin{align}\label{eq:weakFormSi1_F=0}
\int\limits_{t_i}^{t} \langle \partial_t(J_\e(\tau) \hue(\tau),\hue(\tau) \rangle_{\Oe} \dd \tau
+ (A_\e \nabla \hue, \nabla \hue)_{(t_i,t) \times \Oe}
+ (B_\e \hue, \nabla \hue)_{(t_i,t) \times \Oe} = 0,
\end{align}
We note that the left-hand side of \eqref{eq:weakFormSi1_F=0} can be rewritten to
\begin{align}\label{eq:weakFormSi1_F=0_Left}
\int\limits_{t_i}^t \langle \partial_t(J_\e(\tau) \hue(\tau)), \hue(\tau) \rangle_\Oe \dd \tau = \tfrac{1}{2} \norm{\sqrt{J_\e(t)} u\hue(t)}{\Oe}^2 +
\tfrac{1}{2} (\partial_t J_\e \hue, \hue)_{(t_i,t) \times \Oe},
\end{align}
thus \eqref{eq:weakFormSi1_F=0} becomes
\begin{align}\notag
\tfrac{1}{2} \norm{\sqrt{J_\e(t)} \hue(t)}{\Oe}^2 + (A_\e \nabla \hue, \nabla \hue)_{(t_i,t) \times \Oe} 
\\\label{eq:weakFormSi1_F=0'}
=
-(B_\e \hue, \nabla \hue)_{(t_i,t) \times \Oe} - \tfrac{1}{2} (\partial_t J_\e \hue, \hue)_{(t_i,t) \times \Oe}.
\end{align}
Using the uniform boundedness from below of $J_\e$ and the coercivity of $A_\e$ given by Lemma \ref{lemma:Estimates-psi-eps}, we can estimate the left-hand side of \eqref{eq:weakFormSi1_F=0'} by
\begin{align}\label{eq:weakFormSi1_F=0_Est_left}
\tfrac{1}{2}c_J \norm{\hue(t)}{\Oe}^2 +\alpha \norm{\nabla \hue}{(t_i,t)\times\Oe}^2 \leq  \tfrac{1}{2} \norm{\sqrt{J_\e(t)} \hue(t)}{\Oe}^2 + (A_\e \nabla \hue, \nabla \hue)_{(t_i,t) \times \Oe}.
\end{align}
The right-hand side of \eqref{eq:weakFormSi1_F=0'} can be estimated with the Cauchy--Schwarz and Young inequalities for arbitrary $\delta >0$ and a constant $C_\delta$ by
\begin{align}\notag
-(B_\e \hue, \nabla \hue
)_{(t_i,t) \times \Oe} - \tfrac{1}{2} (\partial_t J_\e \hue, \hue)_{(t_i,t) \times \Oe}
\\\label{eq:weakFormSi1_F=0_Est_right}
\leq
\delta \norm{\nabla \hue }{(t_i,t) \times \Omega}^2 + C_\delta \norm{\hue }{(t_i,t) \times \Omega}^2
+
C \norm{\hue }{(t_i,t) \times \Omega}^2
\end{align}
After combining \eqref{eq:weakFormSi1_F=0'}--\eqref{eq:weakFormSi1_F=0_Est_right} and collecting all the constants, we get for $\delta = \alpha /2$
\begin{align*}
\frac{1}{2}c_J \norm{u_\e(t)}{\Oe} + (\alpha-\alpha/2) \norm{\nabla \hue }{(t_i,t) \times \Oe}^2 \leq (C_{\alpha/2} +C)\norm{\hue }{(t_i,t) \times \Oe}^2
\end{align*}
Then, the lemma of Gronwall shows $\hue = 0$ which gives the uniqueness of $\hue$ and thus $L_{\e,2}$ is well defined.

\paragraph{$\bullet$ Uniform bound of $L_{\e,2}$}
In order to derive a uniform bound for $\hue$, we test \eqref{eq:FixedPointUe} with $\chi_{(t_i,t)}\hue$ for a.e.~$t \in S_i$, which gives
\begin{align}\notag
\int\limits_{t_i}^{t} \langle \partial_t(J_\e(\tau) \hue(\tau),\hue(\tau) \rangle_{\Oe} \dd \tau
+ (A_\e \nabla \hue, \nabla \hue)_{(t_i,t) \times \Oe}
+ (B_\e \hue, \nabla \hue)_{(t_i,t) \times \Oe} 
\\\label{eq:weakFormSi1_Tested_ue}
= 
(J_\e\hat{f}^\p_\e, \hue)_{(t_i,t) \times \Oe}  
-\sum\limits_{k \in I_\e} \left( \tfrac{r_{\e,k}^{n-1}}{ r_0^{n-1}}   \e f_\e (\hue, r_{\e,k}),  \hue\right)_{(t_i,t) \times \Gek}
\end{align}
We rewrite the first term of \eqref{eq:weakFormSi1_Tested_ue}, similar to \eqref{eq:weakFormSi1_F=0_Left}, by
\begin{align*}
&\int\limits_{t_i}^t \langle \partial_t(J_\e(\tau) \hue(\tau)), \hue(\tau) \rangle_\Oe \dd \tau 
\\
&= \tfrac{1}{2} \norm{\sqrt{J_\e(t)} \hue(t)}{\Oe}^2 - \tfrac{1}{2} \norm{\sqrt{J_\e(t_i)} \hat{u}_\e^{(t_i)}}{\Oe}^2+
\tfrac{1}{2} (\partial_t J_\e \hue, \hue)_{(t_i,t) \times \Oe}.
\end{align*}
Thus, \eqref{eq:weakFormSi1_Tested_ue} can be rewritten into
\begin{align}\notag
&\tfrac{1}{2} \norm{\sqrt{J_\e(t)} \hue(t)}{\Oe}
+
(A_\e \nabla \hue, \nabla \hue)_{(t_i,t) \times \Oe} 
= 
(J_\e\hat{f}^\p_\e, \hue)_{(t_i,t) \times \Oe}  
\\\notag
&-\sum\limits_{k \in I_\e} \left( \tfrac{r_{\e,k}^{n-1}}{ r_0^{n-1}}   \e f_\e (\hue, r_{\e,k}),  \hue\right)_{(t_i,t) \times \Gek}
-
(B_\e \hue, \nabla \hue)_{(t_i,t) \times \Oe}
\\\label{eq:Weak_Form_rewritten_For_Boundedness}
&-
\tfrac{1}{2} (\partial_t J_\e \hue, \hue)_{(t_i,t) \times \Oe}
+
\tfrac{1}{2} \norm{\sqrt{J_\e(t_i)} \hat{u}_\e^{(t_i)}}{\Oe}^2.
\end{align}
The first two terms of the right-hand side of \eqref{eq:Weak_Form_rewritten_For_Boundedness} can be estimated with the Cauchy--Schwarz and Young inequalities and the $\e$-scaled trace operator \eqref{eq:Eps-scaled-trace} by
\begin{align}\notag
&(J_\e\hat{f}^\p_\e, \hue)_{(t_i,t) \times \Oe}  
-\sum\limits_{k \in I_\e} \left( \tfrac{r_{\e,k}^{n-1}}{ r_0^{n-1}}   \e f_\e (\zeta, r_{\e,k}),  \hue\right)_{(t_i,t) \times \Gek}
\\\notag
&\leq
C\norm{\hat{f}^\p_\e}{(t_i, t)\times \Oe}^2 +\norm{\hue}{(t_i, t)\times \Oe}^2
+
C_\e \e \norm{f_{\max}}{(t_i, t)\times \Ge}^2 +\e\norm{\hue}{(t_i, t)\times \Ge}^2
\\\notag
&\leq
C + \norm{\hue}{(t_i, t)\times \Oe}^2
+
C f_{\max}^2 \e|S_i||\Ge|
+
C_{\delta} \norm{\hue}{((t_i, t)\times \Oe)}^2 + \delta \e \norm{\nabla \hue}{(t_i, t)\times \Oe}^2
\\\label{eq:Est1}
&\leq
C + C_{\delta} \norm{\hue}{(t_i, t)\times \Oe}^2 + \delta \norm{\nabla \hue}{(t_i, t)\times \Oe}^2
\end{align}
Similarly, we obtain
\begin{align}\label{eq:Est2}
&-
(B_\e \hue, \nabla \hue)_{(t_i,t) \times \Oe} 
\leq
C \norm{\hue}{(t_i, t)\times \Oe}^2 +  C \norm{\nabla \hue}{(t_i, t)\times \Oe}^2,
\\
&-
\tfrac{1}{2} (\partial_t J_\e \hue, \hue)_{(t_i,t) \times \Oe}
\leq
C \norm{\hue}{(t_i, t)\times \Oe}^2,
\\\label{eq:Est3}
&\norm{\sqrt{J_\e(t_i)} \hat{u}_\e^{(t_i)}}{\Oe}^2 \leq C \norm{ \hat{u}_\e^{(t_i)}}{\Oe}^2 \leq C_K.
\end{align}
Combining the estimates
\eqref{eq:weakFormSi1_F=0_Est_left}, \eqref{eq:Est1}--\eqref{eq:Est3} with \eqref{eq:Weak_Form_rewritten_For_Boundedness} yields for $\delta$ small enough and after collecting all the constants
\begin{align}
\norm{\hat{u}_\e(t)}{\Oe}^2 + \norm{\nabla \hue}{(t_i, t)\times \Oe}^2
\leq
C_{K}+ C \norm{\hue}{(t_i, t)\times \Oe}^2.
\end{align}
Then, the Lemma of Gronwall implies
\begin{align}\label{eq:est:hue_Si}
\norm{\hat{u}_\e(t)}{\Oe}^2 + \norm{\nabla \hue}{S_i \times \Oe}^2 \leq  C_{K}
\end{align}
for a.e.~$t\in S_i$. 

\paragraph{$\bullet$ Lipschitz estimate of $L_{\e,2}$}

Let $r_{\e,i} \in V_{r,\e}$ and $\zeta_{i} \in L^2(S;H^1(\Oe))$.
We define $\hat{u}_{\e,i} = L_{\e,2}(\zeta_{i}, r_{\e,i})$ for $i \in \{1,2\}$ as well as $\psi_{\e,i}$ and $\Psi_{\e,i}, J_{\e,i}$ by \eqref{eq:def:psi-eps}--\eqref{eq:def:Psi-eps} for $r_\e = r_{\e,i}$ and $A_{\e,i} \coloneqq J_{\e,i} \Psi_{\e,i}^{-1}\Psi_{\e,i}^{-\top}$,  $B_{\e,i} \coloneqq J_{\e,i} \Psi_{\e,i}^{-1} \partial_t \psi_{\e,i}$.
We test \eqref{eq:FixedPointUe} for $i \in \{1,2\}$ with $\chi_{(t_i,t)}(\hat{u}_{\e,1}- \hat{u}_{\e,2})$ and subtract the corresponding equations:
\begin{align}\notag
&\int\limits_{t_i}^t\langle \partial_t (J_{\e,1} \hat{u}_{\e,1} - J_{\e,2} \hat{u}_{\e,2}), \hat{u}_{\e,1} -\hat{u}_{\e,2}\rangle_{\Oe} \dt 
\\\notag
&+ (A_{\e,1} \nabla \hat{u}_{\e,1} - A_{\e,2} \nabla \hat{u}_{\e,2},  \nabla (\hat{u}_{\e,1} -\hat{u}_{\e,2}))_{(t_i,t)\times \Oe} 
\\\notag
&+ (B_{\e,1} \hat{u}_{\e,1}- B_{\e,2} \hat{u}_{\e,2}, \nabla (\hat{u}_{\e,1} -\hat{u}_{\e,2}))_{(t_i,t)\times \Oe} 
\\\notag
&= (J_{\e,1}f^\p(\cdot_t,\psi_{\e,1}(\cdot_t,\cdot_x))- J_{\e,2}f^\p(\cdot_t,\psi_{\e,2}(\cdot_t,\cdot_x)),\hat{u}_{\e,1} -\hat{u}_{\e,2} )_{(t_i,t)\times \Oe} 
\\\label{eq:u-eps-u-eps}
&-\e\sum\limits_{k \in I_\e} \Big(\tfrac{r_{\e,1,k}^{n-1}}{r_0^{n-1}} f(\zeta_{1},r_{\e,1,k}) -\tfrac{r_{\e,2,k}^{n-1}}{r_0^{n-1}} f(\zeta_{2},r_{\e,2,k}), \hat{u}_{\e,1} -\hat{u}_{\e,2} \Big)_{(t_i,t)\times \Gek} 
\end{align}
Employing $r_{\e,1}(t_i) = r_{\e,2}(t_i)$ and $\hat{u}_{\e,1}(t_i) = \hat{u}_{\e,2}(t_i)$, we can rewrite the first term of \eqref{eq:u-eps-u-eps} into
\begin{align*}
&\int\limits_{t_i}^{t}\langle\partial_t (J_{\e,1}(\tau) \hat{u}_{\e,1}(\tau) - J_{\e,2}(\tau) \hat{u}_{\e,2}(\tau)), \hat{u}_{\e,1}(\tau) -\hat{u}_{\e,2}(\tau)\rangle_{\Oe} \dd \tau
\\
&=
\tfrac{1}{2} \norm{\sqrt{J_{\e,1}}(t) (\hat{u}_{\e,1}(t) -\hat{u}_{\e,2}(t))}{\Oe}^2
+
\tfrac{1}{2} (\partial_t J_{\e,1} (\hat{u}_{\e,1}-\hat{u}_{\e,2}),\hat{u}_{\e,1} -\hat{u}_{\e,2})_{(t_i,t)\times \Oe} 
\\
&+
(\partial_t (J_{\e,1} -J_{\e,2}) \hat{u}_{\e,2},\hat{u}_{\e,1} -\hat{u}_{\e,2})_{(t_i,t)\times \Oe}.
\end{align*}
Thus, we can rewrite \eqref{eq:u-eps-u-eps} by:
\begin{align}\notag
&I_1 +I_2 \coloneqq \tfrac{1}{2} \norm{\sqrt{J_{\e,1}}(t) (\hat{u}_{\e,1}(t) -\hat{u}_{\e,2}(t))}{\Oe}^2  
+ (A_{\e,1} \nabla (\hat{u}_{\e,2} - \hat{u}_{\e,2}),  \nabla (\hat{u}_{\e,1} -\hat{u}_{\e,2}))_{(t_i,t)\times \Oe}
\\\notag
&=
-\tfrac{1}{2} (\partial_t J_{\e,1} (\hat{u}_{\e,1}-\hat{u}_{\e,2}),\hat{u}_{\e,1} -\hat{u}_{\e,2})_{\Oe,t} - (\partial_t (J_{\e,1} -J_{\e,2}) \hat{u}_{\e,2},\hat{u}_{\e,1} -\hat{u}_{\e,2})_{(t_i,t)\times \Oe}
\\\notag
&-((A_{\e,1} - A_{\e,2} )\nabla \hat{u}_{\e,2},  \nabla ( \hat{u}_{\e,1} -\hat{u}_{\e,2}))_{\Oe,t}
- (B_{\e,1} \hat{u}_{\e,1}- B_{\e,2} \hat{u}_{\e,2}, \nabla (\hat{u}_{\e,1} -\hat{u}_{\e,2}))_{(t_i,t)\times \Oe}
\\\notag
&+ (J_{\e,1}f^\p(\cdot_t,\psi_{\e,1}(\cdot_t,\cdot_x))- J_{\e,2}f^\p(\cdot_t,\psi_{\e,2}(\cdot_t,\cdot_x)),\hat{u}_{\e,1} -\hat{u}_{\e,2} )_{(t_i,t)\times \Oe}
\\\notag
&- \e\sum\limits_{k \in I_\e} \Big(\tfrac{r_{\e,1,k}^{n-1}}{r_0^{n-1}} f(\zeta_{1},r_{\e,1,k}) -\tfrac{r_{\e,2,k}^{n-1}}{r_0^{n-1}} f(\zeta_{2},r_{\e,2,k}), \hat{u}_{\e,1} -\hat{u}_{\e,2} \Big)_{(t_i,t)\times \Gek}
\\\label{eq:Def:I-1-I_8}
&\eqqcolon
I_3 +I_4 +I_6 +I_7 +I_8
\end{align}
In the next step, we estimate $I_1,I_2$ from below and $I_3, \dots I_8$ from above:

\textbf{$I_1$, $I_2$:}
Lemma \ref{lemma:Estimates-psi-eps} implies:
\begin{align}\label{eq:est:I1}
&\norm{\sqrt{J_{\e,1}}(t) (\hat{u}_{\e,1}(t) - \hat{u}_{\e,2}(t))}{\Oe}^2  \geq c_J \norm{\hat{u}_{\e,1}(t) -\hat{u}_{\e,2}(t)}{\Oe}^2
\\\label{eq:est:I2}
&(A_{\e,1} \nabla (\hat{u}_{\e,1} - \hat{u}_{\e,2}),  \nabla (\hat{u}_{\e,1} -\hat{u}_{\e,2}))_{(t_i,t)\times \Oe} \geq \alpha \norm{\nabla (\hat{u}_{\e,1} -\hat{u}_{\e,2})}{(t_i,t)\times \Oe}^2.
\end{align}

\textbf{$I_3 +I_4$:}
Application of Lemma \ref{lemma:Estimates-psi-eps}, Lemma \ref{lemma:Lipschitz-psi-eps}, \eqref{eq:est:hue_Si}, the Cauchy--Schwarz and the Young inequalities yields:
\begin{align}\notag
&-\tfrac{1}{2} (\partial_t J_{\e,1} (\hat{u}_{\e,1}-\hat{u}_{\e,2}),\hat{u}_{\e,1} -\hat{u}_{\e,2})_{(t_i,t)\times \Oe} -
(\partial_t (J_{\e,1} -J_{\e,2}) \hat{u}_{\e,2},\hat{u}_{\e,1} -\hat{u}_{\e,2})_{(t_i,t)\times \Oe}
\\\notag
&\leq
C\norm{\hat{u}_{\e,1} -\hat{u}_{\e,2}}{(t_i,t)\times \Oe}^2 
+
C\norm{\partial_t (r_{\e,1} -r_{\e,2})}{(t_i,t)}  \norm{\hat{u}_{\e,2}}{L^\infty((t_i,t);L^2(\Oe))}\norm{\hat{u}_{\e,1} -\hat{u}_{\e,2}}{(t_i,t)\times \Oe}
\\\label{eq:est:I3+4}
&\leq C\norm{\hat{u}_{\e,1} -\hat{u}_{\e,2}}{\Oe,t}^2 + C_K \norm{\partial_t (r_{\e,1} -r_{\e,2})}{(t_i,t)}^2 
\end{align}

\textbf{$I_5$:} We estimate similar as  \eqref{eq:est:I3+4} and use the boundedness of the radii and get for every $\delta >0$ a constant $C_{K,\delta}< \infty$ such that:
\begin{align}\notag
&-((A_{\e,1} - A_{\e,2} )\nabla \hat{u}_{\e,2},  \nabla (\hat{u}_{\e,1} -\hat{u}_{\e,2}))_{ (t_i,t)\times\Oe}
\\\notag
&\leq \norm{r_{\e,1}-r_{\e,2}}{L^\infty((t_i,t))} \norm{\nabla \hat{u}_{\e,2}}{(t_i,t) \times\Oe} \norm{\nabla (\hat{u}_{\e,1} -\hat{u}_{\e,2})}{(t_i,t) \times\Oe}
\\\label{eq:est:I5}
&\leq
C_{K,\delta}\norm{r_{\e,1}-r_{\e,2}}{L^\infty((t_i,t))}^2 +\delta \norm{\nabla (\hat{u}_{\e,1} -\hat{u}_{\e,2)}}{(t_i,t) \times\Oe}^2
\end{align}

\textbf{$I_6$:} Similarly to the estimate of $I_5$, we get:
\begin{align}\notag
&-(B_{\e,1} \hat{u}_{\e,1}-B_{\e,2} \hat{u}_{\e,2}, \nabla (\hat{u}_{\e,1} -\hat{u}_{\e,2}))_{(t_i,t)\times\Oe} 
\\
&\leq\notag (\norm{(B_{\e,1}-B_{\e,2})\hat{u}_{\e,1}}{(t_i,t) \times\Oe} 
+ \norm{B_{\e,2}}{L^\infty(S \times \Oe )} \norm{\hat{u}_{\e,1} -\hat{u}_{\e,2}}{(t_i,t)\times\Oe}) \norm{\nabla (\hat{u}_{\e,1} -\hat{u}_{\e,2})}{(t_i,t)\times\Oe}
\\
&\leq\notag
C_{K,\delta}\norm{r_{\e,1} -r_{\e,2}}{L^\infty((t_i,t))}^2 + C_{K,\delta} \norm{\partial_t (r_{\e,1} - r_{\e,2})}{(t_i,t)}^2
\\\label{eq:est:I6}
&
+C_\delta\norm{\hat{u}_{\e,1} -\hat{u}_{\e,2}}{(t_i,t)\times\Oe}^2
+
\delta \norm{\nabla (\hat{u}_{\e,1} -\hat{u}_{\e,2})}{(t_i,t)\times\Oe}^2
\end{align}

\textbf{$I_7:$} By the same procedure as in the estimate of $I_5$ and employing that $f^\p$ is Lipschitz continuous in each $\e$-scaled cell, we can estimate:
\begin{align}\notag
&(J_{\e,1}f^\p(\cdot_t,\psi_{\e,1}(\cdot_t,\cdot_x))- J_{\e,2}f^\p(\cdot_t,\psi_{\e,2}(\cdot_t,\cdot_x)),\hat{u}_{\e,1} -\hat{u}_{\e,2} )_{(t_i,t)\times\Oe}
\\\notag
&=
((J_{\e,1}- J_{\e,2})f^\p(\cdot_t,\psi_{\e,1}(\cdot_t,\cdot_x)) + J_{\e,2}(f^\p(\cdot_t,\psi_{\e,1}(\cdot_t,\cdot_x))- f^\p(\cdot_t,\psi_{\e,2}(\cdot_t,\cdot_x))),\hat{u}_{\e,1} -\hat{u}_{\e,2} )_{(t_i,t)\times\Oe}
\\\notag
&\leq 
C\norm{r_{\e,1}-r_{\e,2}}{L^\infty((t_i,t))}^2+ 
C\norm{\psi_{\e,1}-\psi_{\e,2}}{L^\infty((t_i,t)\times \Oe)}^2+ C\norm{\hat{u}_{\e,1} -\hat{u}_{\e,2}}{(t_i,t)\times\Oe}^2
\\\label{eq::est:I7}
&\leq
C \norm{r_{\e,1}-r_{\e,2}}{L^\infty((t_i,t))}^2+ C\norm{\hat{u}_{\e,1} -\hat{u}_{\e,2}}{(t_i,t)\times\Oe}^2
\end{align}

$I_8$: Using the Cauchy--Schwarz inequality gives
\begin{align}\notag
&\e\sum\limits_{k \in I_\e} \Big(\tfrac{r_{\e,1,k}^{n-1}}{r_0^{n-1}} f(\zeta_{1},r_{\e,1,k}) -\tfrac{r_{\e,2,k}^{n-1}}{r_0^{n-1}} f(\zeta_{2},r_{\e,2,k}), \hat{u}_{\e,1} -\hat{u}_{\e,2} \Big)_{(t_i,t) \times \Gek}
\\\label{eq:est1}
&\leq
\e\sum\limits_{k \in I_\e} \norm{\tfrac{r_{\e,1,k}^{n-1}}{r_0^{n-1}} f(\zeta_{1},r_{\e,1,k}) -\tfrac{r_{\e,2,k}^{n-1}}{r_0^{n-1}} f(\zeta_{2},r_{\e,2,k})}{(t_i,t) \times \Gek}
\norm{\hat{u}_{\e,1}-\hat{u}_{\e,2}}{(t_i,t) \times \Gek}
\end{align}
We estimate the first factor of the right-hand side of  \eqref{eq:est1} using the Lipschitz continuity of $f$ and the boundedness of $r_{\e,1,k}$ and $r_{\e,2,k}$ 
\begin{align}\notag
&\norm{\tfrac{r_{\e,1,k}^{n-1}}{r_0^{n-1}} f(\zeta_{1},r_{\e,1,k}) -\tfrac{r_{\e,2,k}^{n-1}}{r_0^{n-1}} f(\zeta_{2},r_{\e,2,k})}{(t_i,t) \times \Gek}
\leq
f_{\max}\norm{\tfrac{r_{\e,1,k}^{n-1} - r_{\e,2,k}^{n-1}}{r_0^{n-1}}}{(t_i,t) \times \Gek}
\\\notag
&+
\norm{ \tfrac{r_{\e,2,k}^{n-1}}{r_0^{n-1}} (f(\zeta_{1},r_{\e,1,k})  - f(\zeta_{2},r_{\e,1,k}))}{(t_i,t) \times \Gek}
+
\norm{\tfrac{r_{\e,2,k}^{n-1}}{r_0^{n-1}}f(\zeta_{2},r_{\e,1,k})  - f(\zeta_{2},r_{\e,2,k})}{(t_i,t) \times \Gek}
\\\label{eq:est2}
&\leq C_\e \norm{r_{\e,1,k}- r_{\e,2,k}}{L^\infty((t_i,t))}
+
\sqrt{|S_i|} C_\e \norm{\zeta_{1} - \zeta_{2}}{L^2(S;H^1(\Oe))}.
\end{align}
Thereby we have estimated with the H\"older inequality and the Lipschitz continuity of $f$, 
\begin{align*}
&\norm{f(\zeta_{1},r_{\e,1,k})  - f(\zeta_{2},r_{\e,1,k})}{(t_i,t) \times \Gek}
\\
&\leq
\norm{1}{(t_i,t) \times \Gek}  \norm{  (f(\zeta_{1},r_{\e,1,k})  - f(\zeta_{2},r_{\e,1,k}))^2}{(t_i,t) \times \Gek}
\\
&\leq 
\sqrt{|S_i||\Gek|}  2 f_{\max} 
\norm{  f(\zeta_{1},r_{\e,1,k})  - f(\zeta_{2},r_{\e,1,k})}{(t_i,t) \times \Gek}
\\
&\leq 
\sqrt{|S_i|} C_\e
\norm{\zeta_{1} - \zeta_{2}}{(t_i,t) \times \Gek}
\leq
\sqrt{|S_i|} C_\e \norm{\zeta_{1} - \zeta_{2}}{L^2(S;H^1(\Oe))}.
\end{align*}
Combining \eqref{eq:est1}--\eqref{eq:est2} and applying the Young and the trace inequalities yields:
\begin{align}\notag
&\e\sum\limits_{k \in I_\e} \Big(\tfrac{r_{\e,1}^{n-1}}{r_0^{n-1}} f(\zeta_{1},r_{\e,1,k}) -\tfrac{r_{\e,2,k}^{n-1}}{r_0^{n-1}} f(\zeta_{2},r_{\e,2,k}), \hat{u}_{\e,1} -\hat{u}_{\e,2} \Big)_{(t_i,t) \times \Gek}
\\\notag
&\leq
\sum\limits_{k \in I_\e} C_\e \norm{r_{\e,1,k}- r_{\e,2,k}}{L^\infty((t_i,t))}^2
+
|S_i|C_\e \norm{\zeta_{1} - \zeta_{2}}{L^2(S;H^1(\Oe))}^2
\\\label{eq:est:I8}
&+
C_{\e,\delta}\norm{\hat{u}_{\e,1}-\hat{u}_{\e,2}}{(t_i,t) \times \Oe}^2
+
\delta \norm{\nabla\hat{u}_{\e,1}-\hat{u}_{\e,2}}{(t_i,t) \times \Oe}^2
\end{align}
Now, we combine \eqref{eq:Def:I-1-I_8} with \eqref{eq:est:I1}--\eqref{eq:est:I8} and get for $\delta$ small enough after collecting the constants
\begin{align}\notag
&\norm{\hat{u}_{\e,1}(t) -\hat{u}_{\e,2}(t)}{\Oe}^2
+ \norm{\nabla (\hat{u}_{\e,1} -\hat{u}_{\e,2})}{(t_i,t)\times \Oe}^2
\\\notag
&\leq 
C_\e \norm{(\hat{u}_{\e,1} -\hat{u}_{\e,2})}{(t_i,t)\times \Oe}^2
+C_K 
\norm{\partial_t (r_{\e,1} -r_{\e,2})}{(t_i,t)\times \Oe}^2
\\\notag
&+
C_{K,\e} 
\norm{r_{\e,1} -r_{\e,2}}{L^\infty((t_i,t))\times \Oe}^2
+
|S_i|C_\e \norm{\zeta_{1} - \zeta_{2}}{L^2(S;H^1(\Oe))}^2.
\end{align}
Then, the lemma of Gronwall gives
\begin{align}\notag
&\norm{\hat{u}_{\e,1}(t) -\hat{u}_{\e,2}(t)}{\Oe}^2
+ \norm{\nabla (\hat{u}_{\e,1} -\hat{u}_{\e,2})}{(t_i,t)\times \Oe}^2
\\\notag
&\leq 
C_{\e,K} 
\norm{\partial_t (r_{\e,1} -r_{\e,2})}{(t_i,t)\times \Oe}^2
+
C_{K,\e} 
\norm{r_{\e,1} -r_{\e,2}}{L^\infty((t_i,t))\times \Oe}^2
\\\label{eq:LipEstL2}
&+
|S_i|C_\e \norm{\zeta_{1} - \zeta_{2}}{L^2(S;H^1(\Oe))}^2.
\end{align}

\paragraph{$\bullet$ Lipschitz estimate of $L_\e$}
We combine \eqref{eq:LipEstL1a},\eqref{eq:LipEstL1b},\eqref{eq:LipEstL2} and get for $r_{\e,i} \coloneqq L_{\e,1}(\zeta_{i})$ for $i \in \{1,2\}$:
\begin{align*}
&\norm{L_\e(\zeta_{1}) -L_\e(\zeta_{2})}{L^2(S;H^1(\Oe))}^2
=
\norm{L_{\e,2}(\zeta_{1}, r_{\e,1}) -
L_{\e,2}(\zeta_{2}, r_{\e,2})}{L^2(S;H^1(\Oe))}^2
\\
&\leq
C_{K,\e} 
\norm{\partial_t (r_{\e,1} -r_{\e,2})}{(t_i,t)\times \Oe}^2
+
C_{K,\e} 
\norm{r_{\e,1} -r_{\e,2}}{L^\infty((t_i,t))\times \Oe}^2
\\
&+
|S_i|C_\e \norm{\zeta_{1} - \zeta_{2} }{L^2(S;H^1(\Oe))}^2
\leq
|S_i|C_{K,\e} \norm{\zeta_{1} - \zeta_{2}}{L^2(S;H^1(\Oe))}^2 
\end{align*}
Thus, $L_\e$ becomes a contraction for $\sigma_{\e,K} = (4C_{K,\e})^{-1}$ and hence there exists a unique solution of \eqref{eq:def:psi-eps}--\eqref{eq:def:Psi-eps}, 
\eqref{eq:weakFormSi1}--\eqref{eq:weakFormSi2}.
\end{proof}

Rescaling the trace inequality of the reference cell onto $\Oe$ yields for every $\delta >0$ a constant $C_\delta$ such that for every $\e$ and $u\in H^1(\Oe)$
\begin{align}\label{eq:Eps-scaled-trace}
\norm{u}{\partial \Oe}^2
&\leq
\e \delta \norm{\nabla u}{\Oe}^2 + \e^{-1} C_\delta \norm{u}{L^2(\Oe)}^2.
\end{align} 

\section{Derivation of the limit problem for the periodic substitute problem}

We use the notion of two-scale convergence which was introduced in \cite{All92} and \cite{Ngu89}.
\begin{defi}[Two-scale convergence]
Let $p,q,p_s,q_s \in (1,\infty)$ with $\tfrac{1}{p}+\tfrac{1}{q}= 1$ and $\tfrac{1}{p_s}+\tfrac{1}{q_s}= 1$. We say that a sequence $u_\e$ in $L^{p_s}(S;L^p(\Omega))$ two-scale converges weakly to $u_0 \in L^{p_s}(S;L^p(\Omega \times Y))$ if
\begin{align}
\limEps \intSO u_\e\tx \varphi\txxeps dxdt = \intSOY u_0\txy \varphi\txy dydxdt
\end{align}
for every $\varphi \in L^{q_s}(S;L^{q}(\Omega;C_\#(Y)))$. In this case, we write $u_\e \TscW{p_s,p} u_0$.

Moreover, we say that $u_\e$ two-scale converges strongly to $u_0$ if additionally $\limEps\norm{u_\e}{L^{p_s}(S;L^p(\Omega))} = \norm{u}{L^{p_s}(S;L^p(\Omega \times S))}$.  In this case, we write $u_\e \TscS{p_s,p} u_0$.
\end{defi}

The notion of two-scale convergence provides the following compactness results. Proposition \ref{prop:TwoScaleCompactnessLp} and Proposition \ref{prop:TwoScaleCompactnessW1p} are time dependent versions of compactness results that can be found in \cite{All92}.
\begin{prop}\label{prop:TwoScaleCompactnessLp}
Let $p_s,p \in (1,\infty)$ and let $u_\e$ be a bounded sequence in $L^{p_s}(S;L^p(\Omega))$. Then, there exists a subsequence $\e$ and $u_0 \in L^{p_s}(S;L^p(\Omega \times Y))$ such that
$u_\e \TscW{p_s,p} u_0$.
\end{prop}
For the sake of simplicity, let in the following Proposition the domain $\Oe$ be given as in the previous sections and $\Yp = Y\setminus B_{r_0} (x_M)$ (for more general domains cf.~\cite{All92}).
We use $\widetilde{\cdot}$ in order to denote the extension of functions which are defined on $\Oe$ or $\Oe(t)$ by $0$ to $\Omega$. We use it also for the extension by $0$ to $Y$ for functions which are defined on $\Yp$ or on $Y^*_r$ with $r \in[r_{\min},r_{\max}]$.

\begin{prop}\label{prop:TwoScaleCompactnessW1p}	Let $p_s,p \in (1,\infty)$ and let $u_\e$ be a bounded sequence in $L^{p_s}(S;W^{1,p}(\Oe))$. Then, there exists a subsequence $\e$ and $(u_0,u_1) \in L^{p_s}(S;W^{1,p}(\Omega)) \times L^{p_s}(S;L^p(\Omega;W^{1,p}_\#(\Yp)/\R))$ such that
$\widetilde{u_\e} \TscW{p_s,p} \xYp u_0$ and $\widetilde{\nabla u_\e} \TscW{p_s,p} \xYp \nabla_x u_0 + \widetilde{\nabla_y u_1}$.
\end{prop}

In order to have \eqref{diag} commutative, we use  concept of locally periodic transformations, which was introduced for the stationary case in \cite{Wie21} and is extend to the time-dependent case here:
\begin{defi}\label{def:psi-eps}
We say a sequence of $\psi_\e : S \times \Omega \rightarrow \Omega$, is a sequence of locally periodic transformations if
\begin{enumerate}
\item \label{item:def:psi-eps:1} $\psi_\e \in L^\infty(S;C^1(\Omega))^N$
\item \label{item:def:psi-eps:2} there exists a constant $c_J$ such that $J_\e(t) \geq c_J$ for a.e.~$t\in S$ with $J_\e\tx \coloneqq \det(\Psi_\e \tx)$ and $\Psi_\e \coloneqq D_x \psi_\e \tx$,
\item \label{item:def:psi-eps:3} there exists a constant $C >0$ such that $\e^{i-1} \norm{\check{\psi}_\e}{L^\infty(S;C^i(\Omega))} \leq C$ for $i \in \{0,1\}$, where $\check{\psi}_\e\tx \coloneqq \psi_\e\tx -x$ is the corresponding displacement mapping,
\item \label{item:def:psi-eps:4} there exists $\psi_0 \in L^\infty(S \times \Omega;C^1(Y))^N$, which we call limit transformation, such that
\begin{enumerate}
\item \label{subitem:def:psi-eps:1} $\psi_0(t,x,\cdot_y) : Y \rightarrow Y$ are $C^1$-diffeomorphisms for a.e.~$\tx \in \SO$ with inverses $\psi_0^{-1}(t,x,\cdot_y)$ for $\psi_{\e}^{-1} \in L^\infty(\SO; C^1(Y))$,
\item \label{subitem:def:psi-eps:2} the corresponding displacement mapping, defined for a.e.~$\tx \in \SO$ by $\check{\psi}_0\txy \coloneqq \psi_0\txy -y$, can be extended $Y$-periodically such that $\check{\psi}_0 \in L^\infty(\SO;C^1_\#(\overline{Y}))^N$,
\item \label{subitem:def:psi-eps:3} $\e^{-1} \check{\psi}_\e \TscS{p,p} \check{\psi}_0$ and $\nabla \check{\psi}_\e \TscS{p,p} \nabla_y \check{\psi}_0$ for every $p\in (1,\infty)$.
\end{enumerate}
\end{enumerate}
For a.e.~$\tx \in \SO$, we denote the Jacobian matrix and determinant of $y \mapsto \psi_0(t,x,y)$ by $\Psi_0\txy\coloneqq D_y \psi_0 \txy$ and $J_0 \txy \coloneqq \det(\Psi_0 \txy)$. Moreover, we denote the displacement mappings of the back-transformations by $\check{\psi}_\e^{-1}\tx \coloneqq \psi_\e^{-1}\tx-x$ and $\check{\psi}_0^{-1}\txy \coloneqq \psi_0^{-1}\txy -y$.
\end{defi}
\begin{nota}
For a function $u$, we introduce the following notations:
\begin{align*}
&u_{\psi_\e}\tx \coloneqq u(t,\psi_\e\tx), && u_{\psi_\e^{-1}}\tx \coloneqq u(t,\psi_\e^{-1}\tx), \\
&u_{\psi_0}\txy \coloneqq u(t,x,\psi_0\txy), && u_{\psi_0^{-1}}\txy \coloneqq u(t,x,\psi_0^{-1}\txy), 
\end{align*}
\end{nota}
\begin{rmk}
Let $\psi_\e$ be a sequence of locally periodic transformations in the sense of Definition \ref{def:psi-eps}. If additionally $\partial_t \psi_\e \in L^{p_s}(S;C(\Omega))$ for $p_s>1$, we can conclude $\psi_\e \in C(\overline{S};C(\Oe))$, which allows us to evaluate $\Oe(t) = \psi_\e(t, \Oe)$ for every $t \in \overline{S}$.
Moreover, if $\partial_t \psi_0 \in L^\infty(\Omega;L^{p_s}(S;C(Y)))$, we get $\psi \in L^\infty(\Omega;C(S;C(Y)))$ which allows defining $\Ypx(t)$ for a.e.~$x \in \Omega$ and every $t\in \overline{S}$.
\end{rmk}
In our case, where $\psi_\e$ is given by \eqref{eq:def:psi-eps}, we can show that $\psi_\e$ is a locally periodic transformation in the sense of Definition \ref{def:psi-eps} if $r_\e$ converges strongly.
In order to prove this, we use the unfolding operator
\begin{align*}
L^{p_s}(S;L^p(\Omega)) \to L^{p_s}(S;L^p(\Omega \times Y)), \ \Te u \txy \coloneqq u(t, [x]_{\e, Y},\e y).
\end{align*}
It allows us to rewrite, two-scale convergence as convergence in $L^{p_s}(S_i;L^p(\Omega  \times Y)$, i.e.~
$u_\e \TscW{p_s,p} u_0$ if and if only $\Te u_\e \rightharpoonup u_0$ in $L^{p_s}(S;L^p(\Omega \times Y))$. In our case $\Te$ is isometric, because $\Omega$ consist only on whole $\e$-scaled cells. Thus, $\Te$ is isometric and $u_\e \TscS{p_s,p} u_0$
if and only if
$\Te u_\e \rightarrow u_0$ in $L^{p_s}(S;L^p(\Omega \times Y))$. Moreover, the unfolding operator can be defined for the periodic boundary in the same way, i.e.~
\begin{align*}
L^{p_s}(S;L^p(\Ge)) \to L^{p_s}(S;L^p(\Omega \times \Gamma)), \ \Te u \txy \coloneqq u(t, [x]_{\e, Y},\e y).
\end{align*}
In the limit process, we use the following properties of $\Te$ which can be found in \cite{CDD+12}:
For $L^{p_s}(S;W^{1,p}(\Omega))$ it holds $\e^{-1}\nabla_y \Te u = \Te\nabla_x u$ and for $L^{p_s}(S;L^p(\Ge))$ it holds $\e \intS \int\limits_{\Ge} u\tx \dxt = \intS \intO  \intG\Te u \txy \sigma_y \dxt$.

\begin{lemma}\label{lem:psi-locally-periodic}
Let $\psi$ be defined by \eqref{eq:def:psi} and $\psi_\e$ by \eqref{eq:def:psi-eps}, where $R$ fulfils the assumptions \eqref{eq:cond:R1}--\eqref{eq:cond:R4}. Let $r_\mathrm{min}\leq r_{\e,k}(t) \leq r_\mathrm{max}$ for every $k \in I_\e$ and a.e.~$t \in S$ and assume that $r_{\e,k(\cdot_x)}$ converges strongly to $r$ in $L^1(S\times \Omega)$. Then, $\psi_\e$ is a sequence of locally periodic transformations in the sense of Definition \ref{def:psi-eps} with limit transformation
\begin{align}
\psi_0 \txy = \psi(r\tx, y).
\end{align}
\end{lemma}
\begin{proof}
The properties \ref{item:def:psi-eps:1}--\ref{item:def:psi-eps:3} of Definition \ref{def:psi-eps} follow directly from the construction of $\psi_\e$ and the uniform boundedness of $r_{\e,k}(t)$ from above and below. 
The strong convergence of $r_{\e,k}(\cdot_x)$ to $r$ transfers this boundedness to its limit: $r_\mathrm{min}\leq r\tx \leq r_\mathrm{max}$ for a.e.~$\tx \in S\times \Omega$. Thus, the properties \ref{subitem:def:psi-eps:1}--\ref{subitem:def:psi-eps:2} of Definition \ref{def:psi-eps} follows from the construction of $\psi_0$.

It remains to show Property \ref{subitem:def:psi-eps:3}, which is equivalent to the strong convergences $\e^{-1}\Te \check{\psi}_\e \to \check{\psi}_0$ and $\Te \nabla_x \check{\psi}_\e = \e^{-1} \nabla_y \Te \check{\psi}_\e  \to \nabla_y \check{\psi}_0$ in $L^p(\Omega \times Y)$. 
Due to the strong convergence of $r_{\e,k}(\cdot_x)$, we can pass to a subsequence $\e$ such that $r_{\e,k(x)}(t) \to r\tx$ for a.e.~$\tx \in S\times \Omega$. Employing the continuity $r \mapsto \check{\psi}(r,y)$ and $r \mapsto \nabla_y \check \psi(r,y)$, we obtain:
\begin{align*}
&\e^{-1} \Te \check{\psi}_\e \txy = \check{\psi}(r_{\e,k_\e([x]_{\e,Y}  +\e y)}(t), \{[x]_{\e,Y} + \e y \}_{\e,Y}) 
=
\check{\psi}(r_{\e,k_\e(x)}(t), y) \to \check{\psi}(r\tx, y), \\
&\Te \nabla_x \check{\psi}_\e \txy = \e^{-1} \nabla_y\Te \check{\psi}_\e \txy
=
\nabla_y \check{\psi}(r_{\e,k_\e(x)}(t), y) \to \nabla_y \check{\psi}(r\tx, y)
\end{align*}
for a.e.~$\tx \in S\times \Omega$.
Since these functions are uniformly bounded in $L^\infty(S \times \Omega \times Y)$ we obtain the convergence of these functions in $L^p(S \times \Omega\times Y)$. Because, the argumentations holds for every arbitrary subsequence, we obtain the desired convergences.
\end{proof}

We obtain the strong convergence of the Jacobian matrix and determinant of the transformations:
\begin{lemma}\label{lem:Jacobians-Two-scale-Convergence}
Let $\psi_\e$ be a locally periodic transformation in the sense of Definition \ref{def:psi-eps} with limit transformation $\psi_0$ and Jacobian matrices and determinants $\Psi_\e, J_\e,\Psi_0, J_0$. Then, $\Psi_\e \TscS{p,p} \Psi_0$, $J_\e \TscS{p,p} J_0$, $\Psi_\e^{-1} \TscS{p,p} \Psi_0^{-1}$ for every $p \in (1,\infty)$.
\end{lemma}
\begin{proof}
Proposition \ref{prop:Two-scale-trafo-Derivative} is the time-dependent version of \cite[Lemma 3.3]{Wie21} and can be proven analogously to the stationary case there.
\end{proof}
Using the notation of locally periodic transformations, we can show that the two-scale limit and the transformation commutes in the following sense:
\begin{prop}[Two-scale transformation]\label{prop:Two-scale-trafo}
Let $\psi_\e$ be a locally periodic transformation in the sense of Definition \ref{def:psi-eps} with limit transformation $\psi_0$. Let $p_s, p \in (1,\infty)$ and $u_\e,\hat{u}_\e = u_\e(\cdot_t,\psi_\e(\cdot_t,\cdot_x)) \in L^{p_s}(S;L^p(\Omega))$. Then, the following statements hold:
\begin{enumerate}
\item \label{item:Two-scale-trafo:Weak}
$u_\e \TscW{p_s,p} u_0$ for $u_0 \in L^{p_s}(S;L^p(\Omega \times Y))$ if and only if $\hat{u}_\e \TscW{p_s,p} \hat{u}_0$ for $\hat{u}_0 = u_{0,\psi_0}$ and equivalently $u_0 = \hat{u}_{0,\psi_0^{-1}}$,
\item \label{item:Two-scale-trafo:Strong}
$u_\e \TscS{p_s,p} u_0$ for $u_0 \in L^{p_s}(S;L^p(\Omega \times Y))$ if and only if $\hat{u}_\e \TscS{p_s,p} \hat{u}_0$ for $\hat{u}_0 = u_{0,\psi_0}$ and equivalently $u_0 = \hat{u}_{0,\psi_0^{-1}}$.
\end{enumerate}
\end{prop}
\begin{proof}
Statement \ref{item:Two-scale-trafo:Weak} is the time-dependent version of \cite[Theorem 3.8]{Wie21} and statement \ref{item:Two-scale-trafo:Strong} is the time-dependent version of  \cite[Theorem 3.14]{Wie21}. They can be proven analogously to the stationary case there.
\end{proof}

We can apply Proposition \ref{prop:Two-scale-trafo} for functions defined on the porous subset by extending  them by $0$ to $\Omega$. However, this can not be transferred directly to the case of weakly differentiable functions because the extension by $0$ is not regularity preserving. Therefore, we use the following transformation rule for functions defined on the porous domain.

For the sake of simplicity, in the following Proposition let the domains $\Oe$ and $\Oe(t)$ be given as in the previous sections and $\Ypx(t) \coloneqq \psi_0(t,x,\Yp)$.
\begin{prop}[Two-scale transformation of gradients] \label{prop:Two-scale-trafo-Derivative}
Let $\psi_\e$ be a locally periodic transformation in the sense of Definition \ref{def:psi-eps} with limit transformation $\psi_0$. Let $p_s, p \in (1,\infty)$ and $u_\e \in L^{p_s}(S;W^{1,p}(\Oe(t)))$, $\hue = u_\e(\cdot_t,\psi_\e(\cdot_t,\cdot_x) )\in L^{p_s}(S;W^{1,p}(\Oe))$,
where $\Oe(t) \coloneqq \psi_\e(t,\Oe)$ for a.e.~$t\in S$.
Then,
$\widetilde{\nabla u_\e} \TscW{p_s,p} \chi_{Y^*_{(\cdot_x)}(\cdot_t)}(\cdot_y) \nabla_x u_0 + \widetilde{\nabla_y u_1}$ for $(u_0,u_1) \in L^{p_s}(S;W^{1,p}(\Omega)) \times L^{p_s}(S;L^p(\Omega;W^{1,p}_\#(\Ypx(t))/\R))$ if and only if $\widetilde{\nabla \hat{u}_\e} \TscW{p_s,p} \xYp \nabla_x \hat{u}_0 +  \widetilde{\nabla_y\hat{u}_1}$ for $\hat{u}_0 = u_{0}$ and $\hat{u}_1 = u_{1,\psi_0} + \xYp\check{\psi}_0 \cdot \nabla_x u_0$, which is equivalent to $u_1 = \hat{u}_{1,\psi_0^{-1}} + \chi_{Y^*_{(\cdot_x)}(\cdot_t)} \check{\psi}_0^{-1} \nabla_x \hat{u}_0$.
\end{prop}
\begin{proof}
Proposition \ref{prop:Two-scale-trafo-Derivative} is the time-dependent version of \cite[Theorem 3.10]{Wie21} and can be proven analogously to the stationary case there.
\end{proof}

In order to homogenise the non-linear boundary terms of \eqref{eq:WeakUnTrans1}--\eqref{eq:WeakUnTrans2}, we need a strong convergence of $\hue$. This can be achieved by extending the functions with the following result from \cite{Hoe16}:
\begin{prop}\label{prop:Extension}
There exists a family of extension operator $E_\e \in L(H^1(\Oe); H^1(\Omega))$ such that 
\begin{align*}
\norm{E_\e u_\e}{\Omega} \leq C \norm{u_\e}{\Oe},
\ \ \ 
\norm{\nabla E_\e u_\e}{\Omega} \leq C \norm{\nabla u_\e}{\Oe}
\end{align*}
for every $u_\e \in H^1(\Oe)$.
\end{prop}
Applying Proposition \ref{prop:Extension} for a.e.~$t\in S$ gives the following time-dependent version of this extension operator.
\begin{cor}\label{cor:Extension}
Let $p_s \in [1,\infty]$. There exists a family of linear extension operators $E_\e$ from $L^{p_s}(S;H^1(\Oe))$ to $L^{p_s}(S;H^1(\Omega))$  such that
\begin{align}\label{eq:ExtensionLp}
&\norm{E_\e u}{L^{p_s}(S;L^2(\Omega))} \leq C \norm{u}{L^{p_s}(S;L^2(\Oe))},
\\\label{eq:ExtensionH^1}
&\norm{\nabla E_\e u}{L^{p_s}(S;L^2(\Omega))} \leq C \norm{\nabla u}{L^{p_s}(S;L^2(\Oe))},
\\\label{eq:ExtensionL_infty_time}
&\norm{E_\e u(t)}{\Omega} \leq C \norm{u(t)}{\Oe}
\end{align}
for every $u \in L^{p_s}(S;H^1(\Oe))$ and a.e.~$t  \in S$.
\end{cor}

In order to show the strong two-scale convergence of $E_\e \hue$, we show a uniform convergence of $\hue(t+h) -\hue(t)$ to $0$ for $h \to 0$. Then, we can conclude with the compactness result of \cite{Sim86}
We define for time-dependent functions $\varphi$ and $h \in \R$
\begin{align}
\delta_h \varphi = \varphi(t+h) -\varphi(t).
\end{align}
\begin{prop}\label{eq:StrongConv_ue}
Let $E_\e{\hue}$ be the extension of $\hue$, where $\hue$ is given by \ref{thm:Existence-ue}. Then, there exists a subsequence $E_\e{\hue}$ and $\hat{u}_0 \in L^2(S\times \Omega)$ such that this subsequence converges strongly to $\hat{u}_0$ in  $L^2(S\times \Omega)$.
\end{prop}
\begin{proof}
Let $E_\e\hue$ be the extension of $\hue$. By Lemma \ref{lem:UniformConvergence_delta_u_eps} it holds that
\begin{align}
\norm{\delta_h E_\e\hue}{\Omega} \leq C \norm{\delta_h \hue}{\Oe} \to 0
\end{align}
converges uniformly (with respect to $\e$) to zero for $h \to 0$.
Moreover, we can estimate for every $0\leq t_1<t_2 \leq S$ with the H\"older inequality
\begin{align*}
&\norm{\int\limits_{t_1}^{t_2} E_\e\hue(t) \dt}{H^1(\Omega)}^2
= 
\intO \Big(\int\limits_{t_1}^{t_2} E_\e\hue\tx \dt \Big)^2 \dx + \intO \Big(\int\limits_{t_1}^{t_2} \nabla E_\e\hue\tx \dt \Big)^2 \dx
\\
&\leq \intO \norm{1}{S}^2 \int\limits_S (E_\e \hue)^2\tx \dt \dx +\intO \norm{1}{S}^2 \int\limits_S( \nabla E_\e \hue)^2\tx \dt \dx 
\\
&= |S| \norm{E_\e\hue}{L^2(S;H^1(\Omega))}^2 \leq C \norm{\hue}{L^2(S;H^1(\Oe))}^2\leq C.
\end{align*}
Since $\int\limits_{t_1}^{t_2} E_\e\hue(t) \dt$ is uniformly bounded in $H^1(\Omega)$, it is compact in $L^2(\Omega)$.
Thus, we can conclude with of \cite[Theorem 1]{Sim86} that $E_\e\hue$ is compact in $L^2(S;L^2(\Omega)) = L^2(S \times \Omega)$.
\end{proof}

\begin{lemma}\label{lem:UniformConvergence_delta_u_eps}
Let $\hue$ be the solution of \eqref{eq:def:psi-eps}--\eqref{eq:def:Psi-eps}, \eqref{eq:WeakTrans1}--\eqref{eq:WeakTrans3}, then 
\begin{align}
\norm{\delta_h \hue}{(0,T-h) \times \Oe} \to 0
\end{align}
for $h\to 0$ uniformly with respect to $\e$, i.e.~there exists a continuous monotonically decreasing function $\omega :[0,\infty) \to \R$ with $\omega(0) = 0$ such that
\begin{align}
\norm{\delta_h \hue}{(0,T-h) \times \Oe} \leq \omega(h)
\end{align}
for every $\delta >0$.
\end{lemma}
\begin{proof}
First we note that
\begin{align*}
\delta_h (J_\e\hue) =
J_\e \delta_h\hue
+
\delta_hJ_\e \hue(\cdot+h).
\end{align*}
Thus,
\begin{align}\notag
&c_J \norm{\delta_h \hue}{(0,T-h) \times \Oe}^2
\leq
(J_\e\delta_h \hue,\delta_h \hue)_{(0,T-h) \times \Oe}
\\\label{eq:Estimate_delta_h_ue}
&\leq
|(\delta_h (J_\e\hue), \delta_h \hue)_{(0,T-h) \times \Oe}|
+
|(\delta_h J_\e\hue(\cdot +h), \delta_h \hue)_{(0,T-h) \times \Oe}|.
\end{align}
Since $\norm{\partial_t r_\e}{L^\infty(S)}$ Lemma \ref{lemma:Estimates-psi-eps} implies $\norm{\partial_t J_\e}{L^\infty(S_i;C(\Oe))}$ and thus, we can estimate the last term of \eqref{eq:Estimate_delta_h_ue} by
\begin{align*}
|(\delta_h J_\e\hue(\cdot +h), \delta_h \hue)_{(0,T-h) \times \Oe}|
\leq
|Ch(\hue(\cdot +h), \delta_h \hue)_{(0,T-h) \times \Oe}|
\\
\leq Ch \norm{\hue(\cdot_t+h)}{(0,T-h) \times \Oe} 2\norm{\hue}{(0,T-h) \times \Oe}
\leq Ch.
\end{align*}
Hence, it is sufficient to show that $(\delta_h (J_\e\hue), \delta_h \hue)_{(0,T-h) \times  \Oe}$ converges uniformly to zero for $h \to 0$.

We note that we can rewrite the first term in \eqref{eq:WeakTrans1}
for $\varphi \in L^2(S;H^1(\Oe))$ with $\partial_t \varphi \in L^2(S;H^{1}(\Oe))$ by
\begin{align*}
&\intS \langle\partial_t( J_{\e}(t) \hat{u}_{\e}(t)), \varphi \rangle_{\Oe} \dt
\\
&=
-
(\partial_t\varphi ,  J_{\e} \hat{u}_{\e} )_{S \times \Oe} + (J_{\e}(T)\hat{u}_{\e}(T),\varphi(T))_{\Oe}
- (J_{\e}(0)\hat{u}_{\e}(0),\varphi(0))_{\Oe}.
\end{align*}
Now, we assume that $\varphi \in H^1((-h,T);H^1(\Oe))$ with $\varphi(-h) = \varphi(T)= 0$. Then, we test \eqref {eq:WeakTrans1} with $\delta_{-h}\varphi$ and use
\begin{align}
(\partial_t \varphi(\cdot_t-h), J_\e \hue)_{S \times \Oe}
=
(\partial_t \varphi, J_\e(\cdot_t+h) \hue(\cdot_t+h) )_{(-h,T-h)\times \Oe},
\end{align}
which yields
\begin{align}\notag
&( \partial_t\varphi, \delta_h(J_{\e} \hat{u}_{\e}) )_{(0,T-h)\times\Oe}
\\\notag
&
=
-
(\partial_t\varphi,  J_{\e}(\cdot_t+h) \hat{u}_{\e}(\cdot_t+h))_{(-h,0) \times\Oe}
+
(\partial_t\varphi, J_{\e}\hat{u}_{\e})_{(T-h,T)\times\Oe}
\\\notag
&+
(J_{\e}(0)\hat{u}_{\e}(0),\varphi(0))_{\Oe}
- (J_{\e}(T)\hat{u}_{\e}(T),\varphi(T-h))_{\Oe}
\\\notag
&+
(A_{\e}(t) \nabla \hat{u}_{\e}(t), \nabla \delta_{-h}\varphi)_{S\times \Oe} +   (B_{\e}(t)\hat{u}_{\e}(t), \nabla \delta_{-h}\varphi)_{S\times \Oe} 
\\\notag
& - (J_{\e}(t) \hat{f}^\p_\e(t), \delta_{-h}\varphi)_{S\times \Oe}
+
\sum\limits_{k\in I_\e}\tfrac{r_{\e,k}^{n-1}}{r_0^{n-1}}\left(\e f(\hat{u}_{\e}(t),r_{\e,k}(t)),\delta_{-h}\varphi\right)_{S\times \Gek}
\\\label{eq:def_Mi}
&\eqqcolon M_1 + \dots + M_8.
\end{align}

Now we define $\eta(t) = \begin{cases}
0 , & t < 0,\\
1, & 0< t < T \\
0 & t \geq T
\end{cases}$
and choose
\begin{align}\label{eq:phi=delta_u_e}
\varphi(t) = h^{-1} \int\limits_{t}^{t+h} \hue(\tau) \eta(\tau) \dd \tau 
\end{align}
where we implicitly extend $\hue(\tau)$ by $0$ for $\tau >T$ and by $\hue(0)$ for $\tau <0$.
Thus, we get for a.e.~$t\in S$
\begin{align}
\partial_t \varphi(t)= \begin{cases}
h^{-1}\hue(t+h) & t < 0,\\
h^{-1}(\hue(t+h) -\hue(t) )& 0 <t < T-h,\\
h^{-1}\hue(t) & t>T-h.
\end{cases}
\end{align}
Then, the left-hand side of \eqref{eq:def_Mi} can be rewritten by
\begin{align}
( \partial_t\varphi, \delta_h(J_{\e} \hat{u}_{\e}) )_{(0,T-h)\times\Oe} = h^{-1}( \delta_h \hue, \delta_h(J_{\e} \hat{u}_{\e}) )_{(0,T-h)\times\Oe}.
\end{align}
Hence, it is sufficient to show that $ M_1,M_2, \dots, M_8$ are uniformly bounded for $\varphi$ given by \eqref{eq:phi=delta_u_e}.

\paragraph{$\bullet$ $M_1, \dots, M_4$} Since $\norm{\hue}{C^0(\overline{S};L^2(\Oe))} \leq C$, we can estimate
\begin{align}
&M_1 = -
(h^{-1}\hat{u}_{\e}(\cdot_t+h), J_{\e}(\cdot_t+h) \hat{u}_{\e}(\cdot_t+h))_{(-h,0) \times\Oe} \leq C,
\\
&M_2 = (h^{-1}\hue, J_{\e}\hat{u}_{\e})_{(T-h,T)\times\Oe} \leq C,
\\
&M_3 = (J_{\e}(0)\hat{u}_{\e}(0), h^{-1} \int\limits_{0}^{h} \hue(\tau)\dd \tau )_{\Oe} \leq C,
\\
&M_4 =- (J_{\e}(T)\hat{u}_{\e}(T),\int\limits_{T-h}^{T} \hue(\tau)\dd \tau)_{\Oe} \leq C.
\end{align}

\paragraph{$\bullet$ $M_5$, $M_6$ and $M_7$}
We show the estimate for $M_5$. The estimate for $M_6$ follows analogously and the estimate for $M_7$ is similar.
We rewrite 
\begin{align*}
M_5 
=&
h^{-1} \intS (A_{\e}(t) \nabla \hat{u}_{\e}(t),  \int\limits_{t-h}^t \nabla \hue(\tau) \eta(\tau) \dd \tau 
\\&-
h^{-1} \intS (A_{\e}(t) \nabla \hat{u}_{\e}(t),  \int\limits_{t}^{t+h} \nabla \hue(\tau) \eta(\tau) \dd \tau )_\Oe \dt
\eqqcolon M_{5a} + M_{5b}.
\end{align*}
Then, we get with the H\"older inequality 
\begin{align*}
M_{5a} 
&\leq
h^{-1} \int\limits_{0}^h \intS C \norm{\nabla \hat{u}_{\e}(t)}{\Oe} \norm{\nabla \hue(t-h+\tau)\eta(t-h+\tau)}{\Oe} \dt  \dd \tau
\\
&\leq 
Ch^{-1} \int\limits_{0}^h  \norm{\nabla \hat{u}_{\e}}{S \times\Oe} \norm{\nabla \hue(\cdot_t-h+\tau)\eta(\cdot_t-h+\tau) }{S \times\Oe}\dd \tau
\\
&\leq
C \norm{\nabla \hat{u}_{\e}}{S \times\Oe}^2 \leq C
\end{align*}
and by the same argumentation we can estimate $M_{5b}$.

\paragraph{$\bullet$ $I_8$}
We split $I_8$ into two sums as we already did for $I_5$. We show the estimate for the first summand. The estimate for the second summand can be done in the same way. With the $\e$-scaled trace inequality 
\begin{align}
\norm{u}{L^1(\partial \Oe)}
&\leq
C (\norm{\nabla u}{L^1(\Oe)} + \e^{-1} \norm{u}{L^1(\Oe)}),
\end{align}
which can be derived in the same way as \eqref{eq:Eps-scaled-trace}, the uniform bound of $r_{\e}$ and $f$ we get
\begin{align*}
&h^{-1} \sum\limits_{k\in I_\e}\left(\e \tfrac{r_{\e,k}^{n-1}}{r_0^{n-1}} f(\hat{u}_{\e}(t),r_{\e,k}(t)),\int\limits_{t-h}^{t} \hue(\tau) \eta(\tau) \dd \tau\right)_{S\times \Gek}
\\
&\leq
h^{-1} \int\limits_{0}^{h} C \sum\limits_{k\in I_\e}\e \int\limits_{S\times \Gek} |\hue(t-h +\tau,x) \eta(t-h +\tau,x)| \dd \tau \dx \dt
\\
&\leq
C \sum\limits_{k\in I_\e}\e \norm{\hue}{L^1(S\times \Gek)}
\leq
C \e \norm{\hue}{L^1(S\times \Ge)}
\\
&\leq
C \norm{\hue}{L^1(S\times \Oe)}
+
C \e \norm{\nabla \hue}{L^1(S\times \Oe)}
\leq
C \norm{\hue}{S\times \Oe}
+
C \e \norm{\nabla \hue}{S\times \Oe} \leq C.
\end{align*}	
Combining the estimates of $M_1, M_2 \dots, M_8$ shows, that $h^{-1} 
( \delta_h\hat{u}_{\e}, \delta_h(J_{\e} \hat{u}_{\e}) )_{(0,T-h)\times\Oe}$ is uniformly bounded and hence that $( \delta_h\hat{u}_{\e}, \delta_h(J_{\e} \hat{u}_{\e}) )_{(0,T-h)\times\Oe}$ converges uniformly to $0$.
\end{proof}

\begin{thm}\label{thm:Limit-Hat-u}
Let $(\hat{u}_\e,r_\e)$ be the unique solution of \eqref{eq:def:psi-eps}--\eqref{eq:def:Psi-eps}, \eqref{eq:WeakTrans1}--\eqref{eq:WeakTrans3}.
Then, 
\begin{align}\label{eq:TwoScaleConv-hue}
&\widetilde{\hue} \TscS{2,2} \xYp \hu_0 &&\textrm{ with respect to the } L^2-\textrm{norm},
\\
&
\widetilde{\nabla \hue} \TscS{2,2} \xYp \nabla_x \hu_0 + \widetilde{\nabla_y \hu_1}  &&\textrm{ with respect to the } L^2-\textrm{norm},
\\\label{eq:TwoScaleConv-re}
&r_{\e,k(\cdot_x)} \rightarrow  r &&\textrm{ in } L^\infty(S;L^p(\Omega)) \textrm{ for every }p \in [1,\infty),
\\\label{eq:TwoScaleConv-dt_re}
&\partial_t r_{\e,k(\cdot_x)} \rightarrow  \partial_t r &&\textrm{ in } L^p(S\times \Omega)) \textrm{ for every }p \in [1,\infty),
\end{align}
where $(\hu_0, \hu_1, r) \in L^2(S;H^1(\Omega)) \times L^2(S;L^2(\Omega;H^1_\#(\Yp)/\R)) \times W^{1,\infty}(S;L^2(\Omega))$ is the unique solution of
\eqref{eq:WeakTwoScaleLimitTrans1}--\eqref{eq:WeakTwoScaleLimitTrans2}.
\end{thm}

Find $(\hu_0,\hu_1, r) \in L^2(S;H^1(\Omega)) \times L^2(S;L^2(\Omega;H^1_\#(\Yp)/\R)) \times W^{1,\infty}(S;L^2(\Omega))$
such that 
\begin{align}\notag
&\intO (1-V_N(r(0,x))) \hu_0^{(0)}(x) \varphi(0,x) \dx
-
\intSO(1-V_N(r\tx)) \hu_0(t,x) \partial_t\varphi\tx \dxt
\\\notag
&+
\intSOYp A_0 \txy (\nabla_x \hu_0(t,x) + \nabla_y \hu_1(t,x,y) \cdot (\nabla_x \varphi \tx + \nabla_y \varphi_1 \txy ) \dyxt
\\\label{eq:WeakTwoScaleLimitTrans1}
&=
\intSO(1-V_N(r\tx)) f^\p\tx \varphi\tx - \partial_t V_N(r\tx) c_s \varphi\tx \dxt
\\\label{eq:WeakTwoScaleLimitTrans2}
&\intS \intO \partial_t r\tx \phi\tx \dxt = \intS \intO \tfrac{1}{c_s} f(\hat{u}_0\tx, r\tx) \phi\tx dxdt
\end{align}
for every $ (\varphi,\varphi_1, \phi) \in H^1(S\times \Omega) \times L^2(S;L^2(\Omega;H^1_\#(\Yp)/\R)) \times L^2(S\times\Omega)$
with initial values $r(0) = r^{(0)}$ and $\hat{u}^{(0)}_0 = u^{(0)}$.
\begin{proof}[Proof of Theorem \ref{thm:Limit-Hat-u}]

Having the uniform estimates \eqref{eq:estimate-u-eps}, we can apply Proposition \ref{prop:TwoScaleCompactnessW1p}, which gives $\hu_0\in L^2(S;H^1(\Omega))$, $\hu_1 \in L^2(S \times \Omega;H^1_\#(\Yp)/\R)$ such that for a subsequence:
\begin{align}\label{eq:Convergence-ue}
\widetilde{\hat{u}_\e} \TscW{2,2} \xYp \hu_0, \ \ \
\widetilde{\nabla u_\e}
\TscW{2,2} \cYp \nabla_x \hu_0 + \widetilde{\nabla_y \hu_1}. 
\end{align}
Proposition \ref{eq:StrongConv_ue} implies (after passing to a further subsequence and identifying the limits)
\begin{align}\label{eq:Convergence-E-ue}
E_\e \hue \to \hat{u}_0 \in L^2(S\times \Omega)
\end{align}
and thus the first convergence of \eqref{eq:Convergence-ue} is strong.
Moreover, this implies
\begin{align}\label{eq:Convergence-Te-ue}
&\Te E_\e\hue \rightarrow  \hat{u}_0  \textrm{ in } L^2(S \times \Omega \times Y)
\end{align}
and we get with $\Te \nabla E_\e \hue = \e^{-1}\nabla_y \Te \hue$, the isometry of $\Te$ and the uniform boundedness of $\norm{\nabla \hue}{S \times \Omega}$:
\begin{align*}
\norm{\nabla_y \Te E_\e \hue}{S \times \Omega \times Y} = \e \norm{\Te \nabla E_\e \hue}{S \times \Omega \times Y} = \e \norm{\nabla E_\e \hue}{S \times \Omega}\leq C \e \norm{\nabla \hue}{S \times \Omega}\rightarrow 0.
\end{align*}
Thus, we can conclude with the trace operator on $\Gamma$
\begin{align}\notag
&\norm{\Te \hue - \hu_0}{S \times \Omega\times \Gamma} = \norm{\Te E_\e \hue - \hu_0}{S \times \Omega\times \Gamma} 
\\\label{eq:ConvergenceGamma}
&\leq C \norm{\Te E_\e \hue - \hu_0}{S \times \Omega\times 
\Yp}  + C \norm{\nabla_y \Te E_\e  \hue - \nabla \hu_0}{S \times \Omega\times \Yp} \to 0.
\end{align}

In order to pass to the limit $\e \to 0$ in the non-linear bulk and boundary terms, we show the strong convergence $r_{\e,k(\cdot_x)} \to r$ at first.
We define $r \in  W^{1,\infty}(S;L^2(\Omega))$ as the unique solution of \eqref{eq:WeakTwoScaleLimitTrans2} with initial value $r(0) = r^{(0)}$ and $\hu_0$ given by \eqref{eq:Convergence-E-ue}. Then, we test \eqref{eq:WeakTrans2} by $\chi_{(0,t)}(r_{\e,k(x)}-r(x))$ for a.e.~$t\in S$, a.e.~$x \in \e k + \e Y$ and every $k \in I_\e$, integrate over $\e k + \e Y$ and sum over $k \in I_\e$:
\begin{align}\notag
&(\partial_t r_{\e,k(\cdot_x)}, r_{\e,k(\cdot_x)} -r)_{\Omega,t}
\\\notag
&= \int\limits_0^t \intO\tfrac{\e^{-N}}{c_sS_{N-1}(r_0) } \int\limits_{\Gamma_{\e,k(x)}} \e f(\hue(\tau,y),r_{\e,k(x)}) d \sigma_y (r_{\e,k(x)}(\tau)-r(\tau,x)) dxd\tau
\\\label{eq:test-r-eps}
&=\left(\tfrac{1}{c_sS_{N-1}(r_0)} f(\Te \hue ,r_{\e,k(\cdot_x)}), r_{\e,k(\cdot_x)} -r\right)_{\Omega\times \Gamma,t}
\end{align}
We test \eqref{eq:WeakTwoScaleLimitTrans2} with $\chi_{(0,t)} (r_{\e,k(\cdot_x)} -r)$  and subtract it from \eqref{eq:test-r-eps}:
\begin{align*}
(\partial_t (r_{\e,k(\cdot_x)} -r), r_{\e,k(\cdot_x)} -r)_{\Omega,t}
=
\tfrac{1}{c_sS_{N-1}(r_0)} (f(\Te \hue, r_{\e,k(\cdot_x)}) -f(\hu_0,r ),r_{\e,k(\cdot_x)} -r)_{\Omega\times \Gamma,t}
\end{align*}
Then, we rewrite the left-hand side and estimate the right-hand side using the Cauchy--Schwarz inequality, the Lipschitz condition \eqref{eq:f-Lipschitz} and the Young inequality:
\begin{align*}
-\tfrac{1}{2}
\norm{r^{(0)}_{\e,k(\cdot_x)} - r^{(0)}}{\Omega}^2
+
\tfrac{1}{2}
\norm{r_{\e,k(\cdot_x)}(t) - r(t)}{\Omega}^2
\leq
C \norm{\Te \hue -u_0}{\Omega\times \Gamma,t}^2 + C\norm{r_{\e,k(\cdot_x)}  -r}{\Omega,t}^2
\end{align*}
We estimate further with the Lemma of Gronwall and pass to the limit using \eqref{eq:ConvergenceGamma} and the strong convergence of the initial values:
\begin{align*}
\norm{ r_{\e,k(\cdot_x)} - r_0}{L^\infty(S;L^2(\Omega))}^2 \leq C \norm{\Te \hue - u_0}{S \times \Omega \times \Gamma}^2 + C \norm{r^{(0)}_{\e,k(\cdot_x)} - r^{(0)}}{\Omega}^2 \to 0.
\end{align*}
Since $r_{\e,k(\cdot_x)}$ and $r_0$ are uniformly bounded in $L^\infty(\SO)$, we get $\norm{ r_{\e,k(\cdot_x)} - r_0}{L^\infty(S;L^p(\Omega))}$ for every $p \in [1,\infty)$. Thus, Lemma \ref{lem:psi-locally-periodic} shows that $\psi_\e$ are locally periodic transformations in the sense of Definition \ref{def:psi-eps} and we can conclude with Lemma \ref{lem:Jacobians-Two-scale-Convergence} the strong two-scale convergence  $J_\e,\Psi_\e, \Psi_\e^{-1}$, which we need in order to pass to the limit $\e \to 0$ in \eqref{eq:WeakTrans1}. 
Moreover, Definition \ref{def:psi-eps}, Proposition \ref{prop:Two-scale-trafo} and \ref{lem:Jacobians-Two-scale-Convergence} can be also formulated for the two-scale convergence without time parameter (cf.~\cite{Wie21}). Thus, we can conclude the strong two-scale convergence also for the initial data, i.e.~$J_\e(0)$ two-scale converges strongly to $J_0(0)$ and $\widetilde{\hue^{(0)}}$ two scale converges to  $\chi_{\Yp_{(\cdot_x)}(\cdot_t)}(\cdot_y)\hu_0^{(0)}(\cdot_x)$ with $\hat{u}_0^{(0)}= u_0^{(0)}(\cdot_x,\psi_0(0,\cdot_x,\cdot_y))$.

The strong convergence $\partial_t r_{\e,k(\cdot_x)} \to \partial_t r$ follows similarly. By testing \eqref{eq:WeakTrans2} and \eqref{eq:WeakTwoScaleLimitTrans2} with $\partial_t (r_{\e,k(\cdot_x)}-r)$ and then subtracting the equations, we can conclude the strong convergence in $L^2(S\times \Omega)$. Subsequently, the boundedness in $L^\infty(S\times \Omega)$ implies the strong convergence in $L^p(S\times \Omega)$ for every $p \in [1,\infty)$. However, we do not need this strong convergence in order to pass to the limit, although the term $B_\e = J_\e\Psi_\e^{-1}\partial_t\psi_\e$ contains the time derivative of $\psi_\e$. The reason is that $\norm{\partial_t \psi_\e}{L^\infty(S\times \Omega)} \leq 
\e C \norm{\partial_t r_{\e,k(\cdot_x)}}{L^\infty(S\times \Omega)}
$ and thus the boundedness of $\norm{\partial_t r_{\e,k(\cdot_x)}}{L^\infty(S\times \Omega)}$ is already sufficient for the limit process.

In order to pass to the limit in \eqref{eq:WeakTrans1}, we test it by $\varphi(\cdot_t,\cdot_x) + \e \varphi_1(\cdot_t,\cdot_x,\tfrac{\cdot_x}{\e})$ for $(\varphi,\varphi_1) \in C^\infty(S;C^\infty(\Omega)) \times D(S;C^\infty(\Omega;C^\infty_\#(Y)))$ with $\varphi(T)= 0$ and integrate the time derivative term by parts:
\begin{align*}
&\intOe J_\e\txy \hue^{(0)}(x) \left(\varphi(0,x) + \e \varphi_1\left(0,x, \tfrac{x}{\e}\right)\right) \dx\\
&-\intSOe J_\e\tx \hue\tx \left(\partial_t  \varphi\tx+
\partial_t  \varphi_1\txxeps\right) \dxt \\
&+ \intSOe A_\e\txy \nabla \hue\tx \cdot \left(\nabla_x \varphi\tx + \e \nabla_x  \varphi_1\txxeps+ \nabla_y \varphi_1\txxeps\right) dxdt
\\
&+ \intSOe B_\e\txy \hue\tx \cdot \left(\nabla_x \varphi\tx + \e \nabla_x  \varphi_1\txxeps + \nabla_y \varphi_1\txxeps\right) dxdt
\\&= 
\intSOe J_\e\tx \hat{f}_\e^\p \tx \left(\varphi\tx + \e \varphi_1\txxeps \right) \dxt
\\
&-\sum\limits_{k \in I_\e} \intS \int\limits_{\Gek}\e \tfrac{r^{n-1}_{\e,k}(t)}{r_0^{n-1}} f(\hue\tx,r_{\e,k}(t)) (\varphi\tx + \e  \varphi_1\tx) \dxt.
\end{align*}
We rewrite the boundary integral with the unfolding operator $\Te$, so that we can pass to the limit $\e \to 0$ using the strong convergences of $\Te u_\e$ and  $r_{\e,k(\cdot_x)}$ and the continuity of $f$:
\begin{align}\notag
&\sum\limits_{k \in I_\e} \intS \int\limits_{\Gek}\e \tfrac{r^{n-1}_{\e,k}(t)}{r_0^{n-1}} f(\hue\tx,r_{\e,k}(t)) (\varphi\tx + \e  \varphi_1\txy) d \sigma_y\dxt
\\\notag
&=
\intS \intO \intG    \tfrac{ r^{n-1}_{\e,k(x)}(t)}{ r_0^{n-1}}  f(\Te \hue\txy, r_{\e,k(x)}(t))
\\\notag&\left(\Te\varphi\tx + \e  \Te\left(\varphi_1\left(\cdot_t,\cdot_x,\tfrac{\cdot_x}{\e}\right)\right)\txy \right)d \sigma_y dx dt
\\\label{eq:testEps}
&\rightarrow
\intS \intO \intG   \tfrac{ r^{n-1}\tx}{ r_0^{n-1}} f (\hu_0\tx, r\tx))  \varphi\tx d \sigma_y dx dt
\end{align}
Using \eqref{eq:WeakTwoScaleLimitTrans2} and $ S_{n-1}(r)=\partial_r V_N(r)$, we can rewrite the right-hand side of \eqref{eq:testEps}:
\begin{align*}
\intS \intO \intG \tfrac{ r^{n-1}\tx}{ r_0^{n-1}} f (\hu_0\tx, r\tx)  \varphi\tx d \sigma_y dx dt = 
\intS \intO \partial_t V_N(r\tx)c_s\varphi\tx dx dt
\end{align*}

Moreover, the uniform boundedness of $\partial_t r_\e$ given by \eqref{eq:estimate-r-eps} implies $\partial_t \psi_\e \to 0$ in $L^\infty(S \times \Omega)$.
Thus, $B_\e \hue$ vanishes in the limit $\e \to 0$ of \eqref{eq:testEps} and we obtain
\begin{align}\notag
&\intOYp J_0(0,x,y) \hu_0^{(0)}(x) \varphi(0,x) dy dx
-
\intSOYp J_0\txy \hu_0^{(0)}(x) \partial_t\varphi\tx \dyxt
\\\notag
&+
\intSOYp A_0 \txy (\nabla_x \hu_0(t,x) + \nabla_y \hu_1(t,x,y)) \cdot (\nabla_x \varphi \tx + \nabla_y \varphi_1 \txy ) \dyxt
\\ \notag
&=
\intSOYp (J_0 \txy f^\p\tx - \partial_t V_N(r\tx) c_s)   \varphi\tx dx dt
\end{align}
which can be rewritten into  \eqref{eq:WeakTwoScaleLimitTrans1}. By a density argument it holds for every  $(\varphi,\varphi_1) \in H^1(S\times\Omega)) \times L^2(S;L^2(\Omega;H^1_\#(\Yp)/\R))$.
Moreover, it can be shown by a fixed-point argument, which is similar to the proof of Theorem \ref{thm:Existence-ue}, that \eqref{eq:WeakTwoScaleLimitTrans1}--\eqref{eq:WeakTwoScaleLimitTrans2} has a unique solution. Thus, the convergence holds for the whole sequence.
\end{proof}

\section{Back-transformation}
Now we transform the two-scale limit problem back from its substitute domain to its actual two-scale domain and obtain the following  transformation-independent weak two-scale formulation.
\begin{thm}[Two-scale limit problem]\label{thm:2ScaleLimit}
Let $(u_\e,r_\e)$ be the solution of \eqref{eq:def:Oet}, \eqref{eq:WeakUnTrans1}--\eqref{eq:WeakUnTrans3}. 
Then, 
\begin{align}\label{eq:Two-Scale-Conv-U_e}
&\widetilde{u_\e}\TscS{2,2} \chi_{Y^*_{r(\cdot_t,\cdot_x)}}(\cdot_y) u_0  &&\textrm{ with respect to the } L^2-\textrm{norm},
\\\label{eq:Two-Scale-Conv-Nabla_U_e}
&\widetilde{\nabla u_\e} \TscS{2,2} \chi_{Y^*_{r(\cdot_t,\cdot_x)}}(\cdot_y) \nabla_x u_0 + \widetilde{\nabla_y u_1}  &&\textrm{ with respect to the } L^2-\textrm{norm}
\end{align}
and the convergences \eqref{eq:TwoScaleConv-re}--\eqref{eq:TwoScaleConv-dt_re} hold,
where $(u_0, u_1, r) \in L^2(S;H^1(\Omega))\times L^2(S;L^2(\Omega;H^1_\#(\Ypx(t))/\R)) \times L^2(S;L^2(\Omega)) $ is the unique solution of the following weak form:

Find $(u_0, u_1, r) \in L^2(S;H^1(\Omega)) \times L^2(S;L^2(\Omega;H^1_\#(\Ypx(t))/\R)) \times W^{1,\infty}(S;L^2(\Omega))$ with $\partial_t ((1 -V_n(r))u_0) \in L^2(S;H^1(\Omega)')$ such that
\begin{align}\notag
&\intS \langle \partial_t ((1- V_N (r(t))) u_0(t)), \varphi(t) \rangle_\Omega \dt
\\\notag
&+ \intS\intOYpx (\nabla_x u_0\tx + \nabla_y u_1\txy) \cdot (\nabla_x \varphi\tx + \nabla_y \varphi_1\txy) dydxdt\\\label{eq:WeakTwoScaleLimit1}
&= 
\intSO(1-V_N(r\tx)) f^\p\tx\varphi\tx  - \partial_t V_N(r\tx)c_s  \varphi\tx \dxt,
\\\label{eq:WeakTwoScaleLimit2}
&\intS \intO \partial_t r\tx \phi\tx dxdt = \intS \intO \tfrac{1}{c_s} f(u_0\tx, r\tx) \phi\tx dxdt
\end{align}
hold
for every $(\varphi,\varphi_1, \phi) \in L^2(S;H^1(\Omega)) \times L^2(S\times \Omega;H^1(\Ypx(t)))\times L^2(S\times\Omega)$
with initial values $r(0) = r^{(0)}$ and $(1- V_N(r)) u_0)(0) = (1- V_N(r^{(0)})) u_0^{(0)}$.
\end{thm}

\begin{proof}
We test \eqref{eq:WeakTwoScaleLimitTrans1} with $(\varphi, \varphi_{1,\psi_0}+\check{\psi}_0 \cdot \nabla_x \varphi)$ for $(\varphi, \varphi_1) \in C^\infty(S;C^\infty(\Omega)) \times C^\infty(S;C^\infty(\Omega;H^1_\#(Y))$ with $\varphi(T)= 0$. Then, we transform the $\Yp$ integral in \eqref{eq:WeakTwoScaleLimit1} with $\psi_0^{-1}(t,x)$ by
\begin{align*}
&\intS\intOYp A_0\txy (\nabla_x \hu_0\tx + \nabla_y \hu_1\txy) 
\\
&\cdot (\nabla_x \varphi\tx + \nabla_y (\varphi_{1,\psi_0}+\check{\psi}_0\txy \cdot \nabla_x \varphi\tx)) dydxdt
\\
=
&\intS\intOYp (\Psi_{0,\psi_0^{-1}}^{-1}\txy  \nabla_x \hu_0\tx + \nabla_y \hu_{1,\psi_0^{-1}}\txy) \\
&\cdot (\Psi_{0,\psi_0^{-1}}^{-1}\txy  \nabla_x \varphi\tx + \nabla_y (\varphi_1(t,x,y)+\check{\psi}_{0,\psi_0^{-1}}\txy \cdot \nabla_x \varphi\tx)) dydxdt
\end{align*}
Using $\Psi_{0,\psi_0^{-1}}^{-1}\txy  = \1  + \nabla_y \check{\psi}_0^{-1}\txy$, we can rewrite
\begin{align*}
&\Psi_{0,\psi_0^{-1}}^{-1}\txy  \nabla_x u_0\tx + \nabla_y \hu_{1,\psi_0^{-1}}\txy
\\
&=
\nabla_x \hu_0\tx + \nabla_y( \hu_{1,\psi_0^{-1}}\txy + \check{\psi}_0^{-1}\txy \cdot \nabla_xu_0\tx)
\\
&= \nabla_x u_0\tx + \nabla_y u_1\txy)
\end{align*}
for a.e.~$\tx \in S \times \Omega$ and a.e.~$y \in \Ypx(t)$ with  $u_0 = \hu_0$ and $u_1 = \hat{u}_{1,\psi_0^{-1}} + \xYpt\check{\psi}_0^{-1} \nabla_x \hat{u}_0$. Note that $\xYpt = \chi_{Y^*_{r(\cdot_t,\cdot_x)}}$.
Using the fact that
\begin{align*}
\check{\psi}_{0,\psi_0^{-1}}\txy &= \check{\psi}_{0}(t,x,\psi_0^{-1}\txy) = \psi_{0}(t,x,\psi_0^{-1}\txy) - \psi_0^{-1}\txy \\&= y - \psi_0^{-1}\txy= -\check{\psi}_0^{-1}\txy
\end{align*}
we get
\begin{align*}
&\Psi_{0,\psi_0^{-1}}^{-1}\txy \nabla_x \varphi\tx + \nabla_y (\varphi_1\txy +\check{\psi}_{0,\psi_0^{-1}}\txy \cdot \nabla_x \varphi\tx)
\\
&= \nabla_x \varphi\tx + \nabla_y (\check{\psi}_0^{-1}\txy \cdot \nabla_x \varphi\tx + \varphi_1\txy +\check{\psi}_{0,\psi_0^{-1}}\txy \cdot \nabla_x \varphi\tx) 
\\
&= \nabla_x \varphi\tx + \nabla_y \varphi_1\txy.
\end{align*}
Thus, we get
\begin{align*}
&\intS\intOYp A_0\txy (\nabla_x \hu_0\tx + \nabla_y \hu_1\txy) 
\\
&\cdot (\nabla_x \varphi\tx + \nabla_y (\varphi_{1,\psi_0}+\check{\psi}_0\txy \cdot \nabla_x \varphi\tx)) dydxdt
\\
&=
\intS\intOYpx (\nabla_x u_0\tx + \nabla_y u_{1}) 
\cdot (\nabla_x \varphi\tx + \nabla_y \varphi_1\txy) dydxdt,
\end{align*}
which allows us to rewrite \eqref{eq:WeakTwoScaleLimitTrans1} into \eqref{eq:WeakTwoScaleLimit1} after integrating the second term of \eqref{eq:WeakTwoScaleLimitTrans1} by parts with respect to time.
By a density argument \eqref{eq:WeakTwoScaleLimit1} holds for all $(\varphi,\varphi_1) \in L^2(S;H^1(\Omega)) \times L^2(S\times \Omega;H^1(\Ypx(t)))$.
The uniqueness of the solution of \eqref{eq:WeakTwoScaleLimit1}--\eqref{eq:WeakTwoScaleLimit2} can be proven by a similar procedure as in the existence proof of Theorem \ref{thm:Existence-ue}.
Thus, the convergence holds for the whole sequence.
The two-scale-convergences \eqref{eq:Two-Scale-Conv-U_e}--\eqref{eq:Two-Scale-Conv-Nabla_U_e} follows from Proposition \ref{prop:Two-scale-trafo} and Proposition \ref{prop:Two-scale-trafo-Derivative}.
\end{proof}
Note that we formulate the initial condition in Theorem \ref{thm:2ScaleLimit} and
Theorem \ref{thm:MainResult}  only for
$(1- V_N(r)) u_0$ 
and not for $u_0$.
The reason is that $1- V_N(r)$ is a priori not regular enough in space in order to transfer the time regularity of $(1- V_N(r)) u_0$ on $u_0$. However, this is not a drawback since $(1- V_N(r)) u_0$ is the actual physically measurable quantity.

\begin{thm}[Homogenised limit problem]\label{thm:MainResult} 
Let $(u_0, r)$ be the unique solution of the two-scale limit problem \eqref{eq:WeakTwoScaleLimit1}--\eqref{eq:WeakTwoScaleLimit2} given by Theorem \ref{thm:2ScaleLimit}.
Then, it is the unique solution of 
\begin{align}\notag
&\intS \langle \partial_t (1- V_N (r(t)) u_0(t)) ,\varphi(t) \rangle_\Omega \dt
+ (A_\mathrm{hom}(r) \nabla_x u_0, \nabla_x \varphi)_{S\times \Omega}
\\\label{eq:homogenised-eq}
&
=((1-V_N(r)) f^\p - \partial_t V_N(r\tx)c_s,  \varphi)_{S\times \Omega}
\end{align}
and \eqref{eq:WeakTwoScaleLimit2}
for every $ (\varphi, \phi) \in L^2(S;H^1(\Omega)) \times L^2(S\times\Omega)$,
where $A_\mathrm{hom}$ is given by
\begin{align}\label{eq:def-Ahom}
(A_\mathrm{hom})_{ij}(r) \coloneqq \int\limits_{Y^*_r} \delta_{ij} + \partial_{y_i} w_j(r;y)
\end{align}
and $w_j(r)$ is the unique solution in $H^1_\#(Y^*_r)/ \R$ such that
\begin{align}\label{eq:CellProblem}
\int\limits_{Y^*_r} (\nabla_y w_j(r;y) + e_j) \cdot \nabla_y \varphi(y) dy = 0.
\end{align}
for every $\varphi \in H^1_\#(Y^*_r)$.
\end{thm}
\begin{proof}
Choosing $\varphi =0$ in \eqref{eq:WeakTwoScaleLimit1} implies $u_1\txy = \sum\limits_{i=1}^N \partial_{x_j} u_0\tx w_j(r\tx,y)$. Inserting this in in \eqref{eq:WeakTwoScaleLimit1} yields \eqref{eq:homogenised-eq} for $A_\mathrm{hom}$ given by \eqref{eq:def-Ahom}.
\end{proof}

In our model the total mass is given by the sum of the mass in the pore space and the mass in the solid space. Thus, the conservation of mass reads $\partial_t (1-V_N(r)) u_0) + \partial_t V_N(r) c_s)=$ density of external sources.
Testing our limit model \eqref{eq:homogenised-eq} with $\varphi \in C^\infty(S)$ yields exactly this
\begin{align}
\partial_t ((1-V_N(r)) u_0) = (1- V_N(r)) f^p - \partial_t V_N(r) c_s.
\end{align}

\section{Acknowledgements}
We would like to thank M.~Gahn for fruitful discussions.

\printbibliography

@InBook{Kro95,
	author = {J. Kropp},
	title = {Performance Criteria for Concrete Durability},
	chapter = {Relations between transport characteristics and durability},
	publisher = {CRC Press},
	year = {1995},
	doi = {https://doi.org/10.1201/9781482271522}
}

@Book{Bie88,
	author = {Thomas A. Bier},
	title = {Karbonatisierung und Realkalisierung von Zementstein und Beton},
	publisher = {Ph.D. Dissertation, University of Karlsruhe},
	year = {1988},
}

@InProceedings{TSR07,
	author = {Alexandre M. Tartakovsky and Paul Meakin and Timothy D Scheibe and Brian D. Wood},
	title = {A smoothed particle hydrodynamics model for reactive transport and mineral precipitation in porous and fractured porous media},
	booktitle = {Water Resources Research},
	year = {2007},
	volume = {43},
	doi = {https://doi.org/10.1029/2005WR004770},
}

@article{VNo08,
	author = {Tycho L. van Noorden},
	title = {Crystal precipitation and dissolution in a porous medium: Effective equations and numerical experiments},
	journal = {Multiscale Model. Simul.},
	year = {2008},
	volume = {7},
	pages = {1220--1236},
	doi = {https://doi.org/10.1137/080722096},
}

@article{VNP10,
	author = {Tycho L. van Noorden and Iuliu Sorin Pop and Anozie Ebigbo and Rainer Helmig},
	title = {An upscaled model for biofilm growth in a thin strip},
	journal = {Water Resources Research},
	year = {2010},
	volume = {46},
	doi = {https://doi.org/10.1029/2009WR008217},
	number = {6},
}

@article{KAP07,
	author = {George E. Kapellos and Terpsichori S. Alexiou and Alkiviades C. Payatakes},
	title = {Hierarchical simulator of biofilm growth and dynamics in granular porous materials},
	journal = {Advances in Water Resources},
	year = {2007},
	volume = {30},
	doi = {https://doi.org/10.1016/j.advwatres.2006.05.030},
	number = {6-7},
	pages = {1648-1667},
}

@article{TZK02,
	author = {Martin Thullner and Josef Zeyer and Wolfgang Kinzelbach},
	title = {Influence of Microbial Growth on Hydraulic Properties of Pore Networks},
	journal = {Transport in Porous Media},
	year = {2002},
	volume = {49},
	doi = {https://doi.org/10.1023/A:1016030112089},
	pages = {99–122},
}

@article{All92,
	author = {Grégoire Allaire},
	title = {Homogenization and two-scale convergence},
	journal = {SIAM J. Math. Anal.},
	year = {1992},
	volume = {23},
	pages = {1482--1518},
}

@article{Ngu89,
	author = {Gabriel Nguetseng},
	title = {A general convergence result for a functional related to the theory of Homogenization },
	journal = {SIAM J. Math. Anal.},
	year = {1989},
	volume = {20},
	pages = {608--623},
}

@article{Pet07,
	author = {Malte Andreas Peter},
	title = {Homogenisation in domains with evolving microstructure},
	journal = {C.~R.~Mecanique},
	year = {2007},
	volume = {335},
	pages = {357--362},
}

@article{PB09a,
	author = {Malte Andreas Peter and Michael Böhm},
	title = {Multiscale Modelling of Chemical Degradation Mechanisms in Porous Media with Evolving Microstructure},
	journal = {Multiscale Model. Simul.},
	year = {2009},
	volume = {7},
	pages = {1643--1668},
}

@article{EM17,
	author = {Michael Eden and Adrian Muntean},
	title = {Homogenization of a fully coupled thermoelasticity problem for a highly heterogeneous medium with a priori known phase transformations},
	journal = {Math. Methods Appl. Sci.},
	year = {2017},
	volume = {40},
	pages = {3955--3972},
}

@article{GNP21,
	author = {Markus Gahn and Maria Neuss-Radu and Iulio Sorin Pop},
	title = {Homogenization of a reaction-diffusion-advection problem in an evolving micro-domain and including nonlinear boundary conditions},
	journal = {J. Differ. Equations},
	year = {2021},
	volume = {289},
	pages = {95--127},
}

@article{Wie21,
	author = {David Wiedemann},
	title = {The two-scale-transformation method},
	journal = {Asymptotic Analysis},
	year = {2022},
	volume = {Pre-press},
	pages = {1-24},
	doi = {https://doi.org/10.3233/ASY-221766},
}

@article{SRF16,
	author = {Raphael Schulz and Nadja Ray and Florian Frank and Hari Shankar Mahato and Peter Knabner},
	title = {Strong solvability up to clogging of an effective
	diffusion–precipitation model in an evolving
	porous medium},
	journal = {European Journal of Applied Mathematics},
	year = {2016},
	volume = {28},
	pages = {179–207},
	doi = {https://doi.org/10.1017/S0956792516000164},
}

@article{MN20,
	author = {Adrian Muntean and Christos Nikolopoulos},
	title = {Colloidal Transport in Locally Periodic Evolving Porous Media - An Upscaling Exercise},
	journal = {SIAM J. Appl. Math.},
	year = {2020},
	volume = {80},
	pages = {448--475},
}

@book{Hoe16,
	author = {Martin H\"opker},
	title = {Extension Operators for Sobolev Spaces on Periodic Domains, Their Applications, and Homogenization of a Phase Field Model for Phase Transitions in Porous Media},
	year = {2016},
	OPTpublisher = {publisher},
}

@article{CDD+12,
	author = {Doina Cioranescu and Alain Damlamian and Patrizia Donato and Georges Griso and Rachad Zaki},
	title = {The periodic unfolding method in domains with holes},
	journal = {SIAM J. Math. Anal.},
	year = {2012},
	volume = {44},
	pages = {718--760},
	doi = {https://doi.org/10.1137/100817942},
}

@Book{Sho97,
author = {Ralph E. Showalter},
	title = {Montone Operators in Banach Space and Nonlinear Partial Differential Equations},
	publisher = {American Mathematical Society},
	year = {1997},
}

@article{Sim86,
	author = {Jacques Simon},
	title = {Compact sets in the space {$L^p(0,T;B$)}},
	journal = {Annali di Matematica Pura ed Applicata},
	year = {1986},
	volume = {146},
	pages = {65–96},
	doi = {https://doi.org/10.1007/BF01762360},
}

@article{KGB+22,
	author = {Mathis Kelm
	 and Stephan Gaerttner and Carina Bringedal
	  and Bernd Flemisch and Peter Knabner and Nadja Ray},
	title = {Comparison study of phase-field and level-set method for three-phase systems including two minerals},
	journal = {Computational Geosciences},
	year = {2022},
	volume = {26},
	pages = {545--570},
	doi = {https://doi.org/10.1007/s10596-022-10142-w},
}

@article{GFP+20,
	author = {Stephan Gaerttner and Peter Frolkovič and Peter Knabner and Nadja Ray},
	title = {Efficiency and Accuracy of Micro‐Macro Models for Mineral Dissolution},
	journal = {Water Resources Research},
	year = {2020},
	volume = {56},
	pages = {1--23},
	doi = {https://doi.org/10.1029/2020WR027585},
}

@article{GFK+22,
	author = {Stephan Gaerttner and Peter Frolkovič and Peter Knabner and Nadja Ray},
	title = {Efficiency of Micro-Macro Models for Reactive Two-Mineral Systems},
	journal = {SIAM Journal on Multiscale Modeling and Simulation},
	year = {2022},
	volume = {206},
	pages = {433--461},
	doi = {https://doi.org/10.1137/20M1380648},
}
\end{document}